\documentclass[1pt]{elsarticle}
\usepackage{geometry}
\usepackage{graphicx,subfigure}
\usepackage[colorlinks,linkcolor=red,citecolor=blue]{hyperref}
\usepackage{amsmath,amssymb,amsxtra,bm}
\usepackage{mathrsfs,amsfonts,amsthm}
\usepackage{appendix}
\usepackage{epstopdf}
\usepackage{booktabs,multirow}
\usepackage{float}
\usepackage{dsfont}
\usepackage{stmaryrd}
\usepackage{cases}
\usepackage[table]{xcolor}

\newtheorem{lemma}{Lemma}

\newtheorem{remark}{Remark}
\newtheorem{example}{Example}

\journal{??}

\newtheorem{sch}{Scheme}[section]
\theoremstyle{plain}

\theoremstyle{remark}

\def\wh{\widehat}
\def\ol{\overline }
\def\wt{\widetilde}
\def\fr{\frac{1}{2}}

\def\d{\mathrm{d}}
\def\diag{\mathrm{diag}}

\def\bn{\mathbf{n}}
\def\bphi{\bm{\phi}}
\def\bx{\mathbf{x}}

\def\bv{\mathbf{v}}
\def\bz{\mathbf{z}}
\def\bu{\mathbf{u}}

\def\bC{\mathbf{C}}
\def\bD{\mathbf{D}}

\def\bI{\mathbf{I}}
\def\bJ{\mathbf{J}}
\def\bH{\mathbf{H}}
\def\bJ{\mathbf{J}}

\def\bT{\mathbf{T}}
\def\bV{\mathbf{V}}

\def\bU{\mathbf{U}}
\def\sech{\mathrm{sech}}

\def\cH{\mathcal{H}}

\newcommand{\ben}{\begin{eqnarray}}
\newcommand{\een}{\end{eqnarray}}
\newcommand{\beq}{\begin{equation}}
\newcommand{\eeq}{\end{equation}}
\newcommand{\bea}{\begin{array}}
	\newcommand{\eea}{\end{array}}
\newcommand{\bef}{\begin{figure}[H]}
	\newcommand{\eef}{\end{figure}}

\numberwithin{equation}{section}

\begin{document}

\begin{frontmatter}

\title{Efficient energy-preserving numerical approximations for the sine-Gordon equation  with Neumann boundary conditions}


\author[mymainaddress1,mythirdaryaddress2]{Qi Hong}
\author[mysecondaryaddress]{Yushun Wang}
\author[mythirdaryaddress,mythirdaryaddress2]{Yuezheng Gong}
\cortext[mycorrespondingauthor]{Corresponding author}
\ead{gongyuezheng@nuaa.edu.cn}
\address[mymainaddress1]{Beijing Computational Science Research Center, Beijing 100193, China}
\address[mysecondaryaddress]{Jiangsu  Key Laboratory for NSLSCS, School of Mathematical Sciences,  Nanjing Normal University, Nanjing 210023, China}
\address[mythirdaryaddress]{College of Science, Nanjing Uninversity of Aeronautics and Astronautics, Nanjing 210016, China}
\address[mythirdaryaddress2]{{Jiangsu Key Laboratory for Numerical Simulation of Large Scale Complex Systems}	}
\begin{abstract}
We present two novel classes of fully discrete energy-preserving algorithms for the sine-Gordon equation subject to Neumann boundary conditions. The cosine pseudo-spectral method is first used to develop structure-preserving spatial discretizations under two different meshes, which result two finite-dimensional Hamiltonian ODE systems. Then we combine the prediction-correction Crank-Nicolson scheme with the projection approach to arrive at fully discrete energy-preserving methods. Alternatively, we introduce a supplementary variable to transform the initial model into a relaxation system, which allows us to construct structure-preserving algorithms more easily. We then discretize the relaxation system directly by using the cosine pseudo-spectral method in space and the prediction-correction Crank-Nicolson scheme in time to derive a new class of energy-preserving schemes. The proposed methods can be solved effectively by the discrete Cosine transform. Some benchmark examples and numerical comparisons are presented to demonstrate the accuracy, efficiency and superiority of the proposed schemes. 
\end{abstract}


\end{frontmatter}
\section{Introduction}
This paper is devoted to designing structure-preserving algorithms for  the undamped two-dimensional  sine-Gordon (SG) equation subject to Neumann boundary conditions, which reads
\begin{align}\label{model-eq}
\begin{cases}
u_{tt} - \Delta u + \phi(x,y) \sin u = 0,\quad \bx =(x, y ) \in \Omega,\\
u(\bx, 0) = \phi_1(\bx),\quad u_t(\bx,0) = \phi_2(\bx),\\
\nabla u \cdot \bn  = 0,\quad (x, y) \in \partial \Omega,
\end{cases}
\end{align}
where $u(\bx, t)$ represents the wave displacement at position $\bx $ and time $t$,
the function $\phi(x,y)$ can be interpreted as a Josephson current density, $\phi_1$ and $\phi_2$ denote initial waveform and velocity, respectively. The SG equation arises in a variety of application areas in 
science and engineering, such as the motion of a rigid pendulum attached to a stretched wire, solid state physics, nonlinear optics, and the stability of fluid motions \cite{josephson1965supercurrents,dodd1982soliton,drazin1989solitons}.

It is important to note  that the system exhibits a canonical Hamiltonian structure. By introducing the velocity $v = u_t$ and  denote $z = (u, v)^T$, we can reformulate the system as
\begin{align}\label{Hamilton-system}
z_t = J \dfrac{\delta \cH }{\delta z}, \quad J = \left(\begin{matrix}
0&\ 1\\
-1&\ 0
\end{matrix} \right),
\end{align}
where $\delta \cH / \delta z$ denotes variational derivatives with respect to $u$ and $v$.
The Hamiltonian system \eqref{Hamilton-system} possesses two important conservation laws, i.e., the conservation of symplecticity
\begin{align}\label{symplectic-property}
\dfrac{\d}{\d t}\omega =0,\quad 
\omega = \int_{\Omega} \left( \d u \wedge \d v \right) \d \bx,
\end{align}
and the conservation of energy
\begin{align}\label{Hamiltonian-energy}
\dfrac{\d }{\d t} \cH(t) =0,\quad
\cH(t)  = \dfrac{1}{2} \int_{\Omega} \big( v^2 + | \nabla u|^2 + 2\phi( 1 - \cos{u}) \big) \d \bx.
\end{align}
Interested readers are referred to \cite{mclachlan1993symplectic,feng2010symplectic,hairer2006geometric} for the details of symplecticity \eqref{symplectic-property}. In the current paper, we mainly focus on the energy conservation property in \eqref{Hamiltonian-energy}.
Due to the Hamiltonian energy conservation property, one would like to retain this property at the discrete level when developing numerical approximation for the model.
Moreover, if a numerical scheme can warrant the energy conservation property, the dynamics of the initial model would be better captured. These are referred to as geometric integrators or structure-preserving algorithms \cite{feng2010symplectic,hairer2006geometric}.

There have been quite a set of papers in the literatures  discussing how to develop numerical method to solve the sine-Gordon equation.  Here, we briefly recall some well-known numerical methods, such as finite difference method \cite{bratsos2007solution,ben1986numerical,sheng2005numerical}, finite element method \cite{argyris1991finite,wang2011convergence},  meshless method \cite{mirzaei2010meshless,dehghan2010numerical}, spectral or pseudo-spectral method \cite{ablowitz1995numerical,asgari2013numerical}, which  is  a classical technique and widely applied to solve some other PDEs \cite{bao2013exponential,bao2014uniform,gong2018linear} in recently years for its high-order accuracy, and so on.
Among these numerical methods, some projects just focus on the algebraic properties, i.e., stability and convergence etc, but ignore its physical properties, such as the energy conservation property. It is worth noting that structure-preserving algorithms always perform a long-term stability as well as the preservation of conservative quantities, for example, the system energy.

Historically, structure-preserving algorithms  for conservative PDEs, especially Hamiltonian systems, have achieved remarkable success in resolving long time dynamics and conservative properties \cite{hairer2006geometric}. Nowadays, a large
number of structure-preserving algorithms have been developed for the sine-Gordon equation. For instance, 
symplectic and multi-symplectic  methods have been investigated  for this model \cite{mclachlan1993symplectic,wang2003high,chen2006symplectic,mclachlan2014high}, etc. 
With the development of structure-preserving algorithms, energy-preserving methods have attracted a lot of attention. Fortunately,  many numerical techniques have been  investigated to construct energy-preserving methods, such as the discrete variational derivative method \cite{furihata1999finite,Furihata2001}, the discrete gradient method \cite{Dahlby2011}, the average vector field (AVF) method \cite{celledoni2012preserving}, the projection method \cite{hairer2006geometric,calvo2006preservation}, time finite element methods \cite{hansbo2001note}, the Hamiltonian boundary value method \cite{brugnano2010hamiltonian}, energy-preserving exponentially-fitted  methods \cite{miyatake2014energy}, the exponential collocation methods \cite{wang2019exponential}, etc. In the meanwhile, energy-preserving algorithms are widely applied on various  Hamiltonian partial differential equations, such as the Schr\"{o}dinger equation \cite{wangTC2019two}, the KdV equation \cite{brugnano2019energy}, Maxwell's equations \cite{kong2019stable}, etc. In particular, for sine-Gordon model, we refer the reader to  \cite{celledoni2012preserving,fei1991two,shi2017energy,gong2019linearly,cai2019structure} for example. However, almost all existing energy-preserving methods for the SG model 
are either fully nonlinear and implicit, or linear-implicit schemes that preserving a modified energy conservation law based on energy quadratization (EQ) technique. These motivate us to develop newly energy-preserving algorithms, which not only reduce the computational cost in numerical simulations but also warrant the original energy conservation law.

In this paper, we propose two new classes of energy-preserving algorithms for the SG equation subject to Neumann boundary conditions. In order to develop spatial high-order structure-preserving algorithms, the cosine pseudo-spectral method is first studied systematically under the half-point grid and the integer grid, respectively, which derive two energy-preserving spatial approximations for the SG equation. The resulting semi-discrete schemes are recast as a canonical Hamiltonian system, which warrants the energy conservation law at the semi-discrete level. Then we combine a prediction-correction Crank-Nicolson method with an energy projection technique in time for the semi-discrete systems to arrive at fully discrete schemes, where the energy conservation law is conserved at the fully discrete level. An alternative approach is to transform the SG model into a relaxation system by the supplementary variable method (SVM), which allows us to enforce the energy conservation property. Furthermore, we discretize the relaxation system by applying the cosine pseudo-spectral method in space and the prediction-correction Crank-Nicolson scheme in time to obtain a new class of energy-preserving schemes. To our surprise,
the proposed methods can be solved efficiently by exploiting the relationship between the spectral differential matrix and the discrete cosine transform (DCT).

In summary, our proposed schemes enjoy the following advantages: (1) the high-order spatial discretization of the SG model with Neumann boundary conditions is still a canonical Hamiltonian system; (2) they can guarantee the original energy conservation law very accurately, as opposed to linear-implicit energy-preserving schemes that preserve a modified energy conservation based on SAV approach \cite{cai2019structure}; (3) compared with these linear-implicit energy-preserving algorithms via EQ or SAV approaches,
they do not need the the nonlinear part of the energy to be bounded from below;  and (4) compared with some traditional structure-preserving algorithms, they greatly reduce the computational cost, where only a scalar nonlinear algebraic equation needs to be solved. In addition, the SVM provide a new paradigm to develop energy-preserving algorithms. Finally, extensive numerical examples involving comparisons with other classical energy-preserving methods are presented to demonstrate the accuracy and efficiency of the proposed algorithms.

The rest of this paper is organized as follows. In Section \ref{sec:SPSD}, the cosine pseudo-spectral method is presented for the SG equation with Neumann boundary conditions, and the relationship between the cosine pseudo-spectral differential matrix and DCT is established. In Sections \ref{sec:TA} and \ref{sect:SVM}, we propose two novel classes of fully discrete energy-preserving methods, including the projection approach and the SVM, respectively. Numerical results are reported in Section \ref{sec:NR}.  We draw some concluding remarks in the end.

\section{Structure-preserving spatial discretization}\label{sec:SPSD}
In this section, we present the cosine pseudo-spectral spatial approximations for the SG model subject to the homogeneous Neumann boundary condition. A rectangular domain $\Omega = [a, b] \times [c, d]$ is uniformly partitioned with mesh size $h_x = (b - a) / N_x$ and $h_y = (d - c) / N_y$, where $N_x$ and $N_y$ are two positive integers.
\subsection{Spatial discretization  on mid-point grids}
In this subsection, we consider the cosine pseudo-spectral method on the following mid-point  grids
\begin{align*}
\Omega^c_h = \bigg\{(x_j,y_k)\big|x_j = a+(j+\frac{1}{2})h_x,\;  y_k = c+(k+\frac{1}{2})h_y,\;  0\leq j\leq N_x - 1,\; 0\leq k\leq N_y-1 \bigg\}.
\end{align*}
Let $V^c_h = \{ u=(u_{j ,k})| (x_j,y_k)\in\Omega^c_h \}$ be the space of cell-center grid functions defined on $\Omega^c_h$. For any two grid functions $u, v \in V^c_h$, we define the following discrete $L^2$ inner product and the corresponding norm
\begin{align}\label{discr-L2-inner}
(u,v)_h = h_xh_y\sum_{j=0}^{N_x-1}\sum_{k=0}^{N_y-1}u_{j,k}v_{j,k},\quad
\|u\|_h = \sqrt{(u,u)_h}.
\end{align}

Define \cite{bao2013exponential,bao2014uniform,gong2014multi,	shen2011spectral}
\begin{align*}
S_N = \{ X_j(x) Y_k(y)| \; j = 0,\cdots,N_x-1, k = 0,\cdots,N_y-1 \}
\end{align*}
as the interpolation space, where $X_j(x)$ and $Y_k(y)$ are the interpolation basis functions given by
\begin{align*}
X_j(x) &= \dfrac{2}{N_x}\sum_{m=0}^{N_x-1}\dfrac{1}{a_m}\cos (m \mu_x (x_j - a) )\cos(m \mu_x (x - a)),\\
Y_k(y) &= \dfrac{2}{N_y}\sum_{m=0}^{N_y-1}\dfrac{1}{b_m}\cos (m \mu_y (y_k - c) )\cos(m \mu_y (y - c)),
\end{align*}
where $\mu_x = \pi/(b-a)$, $\mu_y = \pi/(d-c)$ and
\begin{align*}
a_m =
\begin{cases}
2, \;& m = 0,\\
1,\;  & m \neq 0,
\end{cases}\quad
b_m =
\begin{cases}
2, \; & m = 0,\\
1,\;  & m \neq 0.
\end{cases}
\end{align*}
We define interpolation operator $I_N: C(\Omega) \longrightarrow S_N$ as follows:
\begin{align}\label{Inp-u-eq-NS-NS}
I_Nu(x,y) = \sum_{j=0}^{N_x-1}\sum_{k=0}^{N_y-1}u_{j,k}X_j(x)Y_k(y),
\end{align}
where $u_{j,k} = u(x_j,y_k)$. Furthermore, the second-order cosine pseudo-spectral differentiation matrices $\bD_2^x$ and $\bD_2^y$ can be computed by
\begin{align*}
(\bD_2^x)_{j,m} = \dfrac{\d ^2 X_m(x_j)}{\d x^2},\quad
(\bD_2^y)_{k,m} = \dfrac{\d ^2 Y_m(y_k)}{\d y^2}.
\end{align*}

\begin{lemma}\label{dct-transform}\cite{shen2011spectral}
	For the matrices $\bD^{\alpha}_2 \; (\alpha = x\; \text{or}\; y )$,  there exists the following relation
	\begin{align}\label{dct-eq}
	\bD_2^{\alpha} = \bC_{N_{\alpha}} \mathbf{\Lambda}^{\alpha}_{2} \bC^{-1}_{N_{\alpha}},\quad
	\end{align}
	where $\bC_{N_{\alpha}}$ denotes the discrete cosine transform (DCT-3) with elements
	\begin{align*}
	(\bC_{N_{\alpha}})_{j,m} = \sqrt{\dfrac{2}{N_{\alpha}a_{m} }}\cos \dfrac{m(j+\frac{1}{2})\pi}{N_{\alpha}},\quad
	\end{align*}
	and
	\begin{align*}
	\mathbf{\Lambda}^{\alpha}_{2} = \diag(   \lambda_{\bD^{\alpha}_2,0},  \lambda_{\bD^{\alpha}_2,1}, \cdots, \lambda_{\bD^{\alpha}_2,N_{\alpha} - 1}   ),\quad
	\lambda_{\bD^{\alpha}_2,j } = -(j\mu_{\alpha})^2.
	\end{align*}
	In addition, we have the following relationship
	\begin{align*}
	\bC^T_{N_{\alpha}} = \bC^{-1}_{N_{\alpha}}.
	\end{align*}
\end{lemma}

\begin{remark}
	Due to the relationship \eqref{dct-eq}, we don't need to know the concrete elements of $\bD_2^x$ and $\bD_2^y$ in numerical calculation, since we can evaluate the derivatives by using the DCT-3 algorithm instead of the cosine spectral differentiation matrix.
\end{remark}

Next, we shall present the cosine pseudo-spectral method for the system \eqref{model-eq} at mid-point grid as follows
\begin{align}\label{semi-disc-SG-model-NS-NS}
\begin{cases}
\dfrac{\textrm{d}}{\textrm{d}t} u = v,\\[0.2cm]
\dfrac{\textrm{d}}{\textrm{d}t} v = \Delta_h u - \phi \sin{u},\\
\end{cases}
\end{align}
where $u,v,\phi \in V^c_h$, and $\Delta_h u = \bD_2^x u + u (\bD_2^y)^T \in V^c_h$. 
Let $\bu$, $\bv$ and $\bphi$ be  $N_x \times N_y$ vectors that are arranged in columns by the matrix variable $u$, $v$ and $\phi$, respectively. Then the semi-discrete scheme \eqref{semi-disc-SG-model-NS-NS} can be rewritten a canonical Hamiltonian structure, namely,
\begin{align}\label{semi-discrete-Ham-NS-NS}
\dfrac{\d}{\d t}{\bf z} = \bJ \nabla_{\bz} \bH, \quad
\bJ = \left(\begin{matrix}
{\bf 0}&\ \bI\\
-\bI&\ {\bf 0}
\end{matrix}\right),
\end{align}
where $\bz = (\bu, \bv)^T$, and the Hamiltonian energy $\bH$ is given by
\begin{align}\label{Hamiltonian-NS-NS}
\bH = \dfrac{1}{2}\left( \bv^T \bv - \bu^T \bD \bu + 2 \bphi^T (1 - \cos{\bu}) \right),
\end{align}
with $\bD = \bI_y \otimes \bD_2^x + \bD_2^y \otimes \bI_x$. Here we note that $\bD$ is a symmetric matrix because $\bD_2^x$ and $\bD_2^y$ are symmetric. As we all know, the semi-discrete Hamiltonian system \eqref{semi-discrete-Ham-NS-NS} possesses the energy conservation law
\begin{align}
\dfrac{\d}{\d t} \bH = (\nabla_{\bz}\bH)^T \dfrac{\d}{\d t}\bz = (\nabla_{\bz}\bH)^T \bJ  \nabla_{\bz} \bH =  0.
\end{align}
Note that the discrete energy of the SG equation is defined as
\begin{align}\label{energy_mid}
\cH_h = \frac{1}{2}\|v\|_h^2 - \frac{1}{2}(u,\Delta_h u)_h + (\phi, 1-\cos{u})_h = h_x h_y \bH,
\end{align} which is conservative for the semi-discrete system \eqref{semi-disc-SG-model-NS-NS}.

\subsection{Spatial discretization on regular grids}
In this subsection, we focus on studying the cosine pseudo-spectral method on the following regular grids
\begin{align*}
\Omega^v_h = \bigg\{(x_j,y_k)\big|x_j = a+jh_x,\;  y_k = c+ k h_y,\;  0\leq j\leq N_x ,\; 0\leq k\leq N_y \bigg\}.
\end{align*}
Let $V^v_h = \{ u=(u_{j, k})| (x_j,y_k)\in\Omega^v_h \}$ be the space of vertex grid functions defined on $\Omega^v_h$. 
We define  
\begin{align*}
\overline{S}_N = \mathrm{span}\{ \overline{X}_j(x) \overline{Y}_k(y) | j = 0, \cdots, N_x,\;  k = 0, \cdots, N_y   \}
\end{align*}
as the interpolation space, where $\overline{X}_j(x)$ and $\overline{Y}_k(y)$ are explicitly given by
\begin{align*}
\overline{X}_j(x)  &= \dfrac{2}{N_x}\sum_{m=0}^{N_x}\dfrac{1}{a_m a_j}\cos (m \mu_x (x_j - a) )\cos(m \mu_x (x - a)),\\
\overline{Y}_k(y) &= \dfrac{2}{N_y}\sum_{m=0}^{N_y}\dfrac{1}{b_m b_k}\cos (m \mu_y (y_k - c) )\cos(m \mu_y (y - c)),
\end{align*}
where $\mu_x = \pi/(b-a)$, $\mu_y = \pi/(d-c)$ and
\begin{align*}
a_m = 
\begin{cases}
2,\; m = 0\;\text{or}\;  N_x,\\
1,\; \text{otherwise},
\end{cases}\quad
b_m = 
\begin{cases}
2,\; m = 0\; \text{or}\;  N_y,\\
1,\; \text{otherwise}.
\end{cases}
\end{align*}
Similarly, we can define the interpolation operator $I_N: L^{2}(\Omega) \longrightarrow \overline{S}_N$ as follows:
\begin{align}\label{Inp-u-eq-N-N}
I_Nu(x,y) = \sum_{j=0}^{N_x}\sum_{k=0}^{N_y}u_{j,k}\overline{X}_j(x)\overline{Y}_k(y).
\end{align}
The corresponding second-order cosine pseudo-spectral differentiation matrices are denoted as ${\overline\bD}_2^x$ and ${\overline\bD}_2^y$, which are obtained by 
\begin{align*}
({\overline\bD}_2^x)_{j,m} = \dfrac{\d ^2 \overline{X}_m(x_j)}{\d x^2},\quad
({\overline\bD}_2^y)_{k,m} = \dfrac{\d ^2 \overline{Y}_m(y_k)}{\d y^2}.
\end{align*}
It is not difficult to prove that $\overline{\bD}_2^x$ and $\overline{\bD}_2^y$ are not symmetric matrices. However, we note that they can be transformed into symmetric matrices with a simple transformation.

\begin{lemma}\cite{shen2011spectral}
	For the matrices $\overline{{\bD}}^{\alpha}_2 \; (\alpha = x\; \text{or}\; y )$,  there exists the following relation
	\begin{align}\label{dct-eq2}
	\bT_{\alpha}^{-1}\overline{\bD}_2^{\alpha}\bT_{\alpha} = \bC_{N_{\alpha}} \mathbf{\Lambda}^{\alpha}_{2} \bC^{-1}_{N_{\alpha}}, 
	\end{align}
	where $\bT_{\alpha} = \diag(\sqrt{a_j})_{j=0}^{N_{\alpha}}$ and
	$\bC_{N_{\alpha}}$ denotes the discrete cosine transform (DCT-1) with elements
	\begin{align*}
	(\bC_{N_{\alpha}})_{j,m} = \sqrt{\dfrac{2}{N_{\alpha}a_{j} a_m }}\cos \dfrac{j m \pi}{N_{\alpha}},
	\end{align*}
	and
	\begin{align*}
	\mathbf{\Lambda}^{\alpha}_{2} = \diag(   \lambda_{\overline{\bD}^{\alpha}_2,0},  \lambda_{\overline{\bD}^{\alpha}_2,1}, \cdots, \lambda_{\overline{\bD}^{\alpha}_2,N_{\alpha} }   ),\quad
	\lambda_{\bD^{\alpha}_2,j } = -(j\mu_{\alpha})^2.
	\end{align*}
	Moreover, we have the following relationship
	\begin{align*}
	\bC^T_{N_{\alpha}} = \bC_{N_{\alpha}} = \bC^{-1}_{N_{\alpha}}.
	\end{align*}
\end{lemma}
\begin{remark}
	Similarly,  with the help of  \eqref{dct-eq2}, we can  evaluate the derivatives by using the DCT-1 algorithm instead of the cosine spectral differentiation matrix.
\end{remark}

What follows,  we shall present the cosine pseudo-spectral method for \eqref{model-eq} with homogeneous Neumann boundary  on the regular grid as follows:
\begin{align}\label{semi-disc-SG-model-N-N}
\begin{cases}
\dfrac{\textrm{d}}{\textrm{d}t}  u = v,\\[0.2cm]
\dfrac{\textrm{d}}{\textrm{d}t}  v = \overline{\Delta}_h u - \phi \sin{u},
\end{cases}
\end{align}
where $u,v,\phi \in V^v_h$, and $\overline{\Delta}_h u= \overline\bD_2^x u + u (\overline\bD^y_2)^T \in V^v_h$. Let $U = \bT_x^{-1} u \bT_y^{-1}$, $V = \bT_x^{-1} v \bT_y^{-1}$. Suppose that $\bU$ and $\bV$ are $(N_x+1) \times (N_y+1)$ vectors that are arranged in columns by the matrix variable $u$ and $v$. Multiplying both sides of \eqref{semi-disc-SG-model-N-N} by $\bT_x^{-1}$ and $\bT_y^{-1}$, then we have
\begin{align*}
\begin{cases}
\dfrac{\textrm{d}}{\textrm{d}t} \bU = \bV,\\[0.2cm]
\dfrac{\textrm{d}}{\textrm{d}t} \bV = \overline{\bD} \bU - \nabla_{\bU} F,
\end{cases}
\end{align*}
where $\overline{\bD} = \bI_y \otimes \left( \bC_{N_{x}} \mathbf{\Lambda}^{x}_{2} \bC^{-1}_{N_{x}}  \right) +\left( \bC_{N_{y}} \mathbf{\Lambda}^{y}_{2} \bC^{-1}_{N_{y}}  \right)  \otimes \bI_x$, and
$$F(U) = \sum_{j=0}^{N_x}\sum_{k=0}^{N_y} \left(\dfrac{1}{(T_x)_{j,j} (T_y)_{k,k}}\right)^2 \phi_{j,k} \Big(1-\cos{\big((T_x)_{j,j} (T_y)_{k,k} U_{j,k}\big)}\Big).$$
Therefore, the semi-discrete scheme \eqref{semi-disc-SG-model-N-N} can be written into the following canonical Hamiltonian structure
\begin{align}\label{semi-discrete-Ham-N-N}
\dfrac{\d}{\d t}\left(
\begin{matrix}
\bU \\
\bV
\end{matrix}
\right) = 
\left(\begin{matrix}
{\bf 0}&\ \bI\\
-\bI&\ {\bf 0}
\end{matrix}\right) 
\left(
\begin{matrix}
\nabla_{\bU} \bH \\
\nabla_{\bV} \bH
\end{matrix}
\right),
\end{align}
with the Hamiltonian energy
\begin{align}\label{Hamiltonian-N-N}
\bH = \frac{1}{2} \bV^T \bV - \frac{1}{2} \bU^T \overline{\bD} \bU + F(\bU).
\end{align}
It is obvious to obtain the Hamiltonian energy conservation
\begin{align*}
\bH(t) \equiv  \bH(0).
\end{align*}

\begin{remark}
	For the regular grids, we define a new discrete inner product and the corresponding norm as follows
	\begin{align}\label{new-def-norm}
	\left<u, v \right>_T = h_x h_y\sum_{j=0}^{N_x}\sum_{k=0}^{N_y} \left(\dfrac{1}{(T_x)_{j,j} (T_y)_{k,k}}\right)^2 u_{j,k} v_{j,k}, \quad \|u\|_T = \sqrt{ \left<u, u \right>_T }.
	\end{align}
	We note that this discrete inner product can be derived by using trapezoidal integral formula. Then the discrete energy of the SG equation is defined as 
	\begin{align}\label{energy_regular}
	\cH_h = \frac{1}{2}\|v\|_T^2 - \frac{1}{2} \left<u, \overline{\Delta}_h u \right>_T + \left<\phi, 1-\cos{u} \right>_T.
	\end{align}
	It is readily to show that $\cH_h = h_x h_y \bH,$ and thus the semi-discrete system \eqref{semi-disc-SG-model-N-N} conserves the discrete energy conservation law.	 
\end{remark}

\begin{remark}
	The semi-discrete  systems \eqref{semi-disc-SG-model-NS-NS} and \eqref{semi-disc-SG-model-N-N} on two different meshes can be recast as a canonical Hamiltonian system. Therefore, all existing symplectic integrators and energy-preserving methods can be directly applied to develop the corresponding structure-preserving algorithms, but the resulting fully discrete schemes are often fully nonlinear and implicit. 
	Solving them requires considerable cost. In this paper, we pay attention to developing some new energy-preserving algorithms, which can greatly improve the numerical implementation and reduce the calculation cost.
\end{remark}

\section{Energy-preserving methods based on the projection approach}\label{sec:TA}
In the previous section, the cosine pseudo-spectral method is shown to produce the ODE system with the energy conservation law, which provides a prerequisite for the application of projection approach. In this section, we focus on developing a projection-type energy-preserving method ({\bf PEPM}) for the semi-discrete systems \eqref{semi-disc-SG-model-NS-NS} and \eqref{semi-disc-SG-model-N-N}, which are named {\bf PEPM-M} and {\bf PEPM-R}, respectively. 

For a positive integer $N_t$, we define time step $\tau = T/N_t$, $t_n = n\tau$, $0 \leq n \leq N_t$, and denote
\begin{align*}
\delta_t^+u^n = \dfrac{1}{\tau} \left( u^{n+1} - u^n  \right),\;
u^{n + \fr} = \fr( u^{n+1} + u^n), \;
\wh{u}^{n + \fr} = \fr (\wh{u}^{n + 1} + u^n), \;
\overline{u}^{n+\fr} = \fr (3u^n - u^{n-1}).
\end{align*}
Employing a prediction-correction Crank-Nicolson scheme with the energy projection strategy for the semi-discrete systems \eqref{semi-disc-SG-model-NS-NS} and \eqref{semi-disc-SG-model-N-N}, we obtain the following fully discrete energy-preserving schemes.

\begin{sch}[{\bf PEPM-M}]\label{energy-preserving-projection}
	Given $u^{n-1}, ~u^n, ~v^n \in V_{h}^c, ~ \forall\; n \geq 1,$ we obtain $u^{n+1}, ~v^{n+1} \in V_{h}^c$ through the following two steps:
	\begin{enumerate}
		\item 
		We first compute $\wt{u}^{n+1}$ and $\wt{v}^{n+1}$ via a prediction-correction Crank-Nicolson scheme
		\begin{align}
		\begin{cases}
		\dfrac{\wh{u}^{n + 1} - u^n}{\tau} = \wh{v}^{n+\fr},\\[0.3cm]
		\dfrac{\wh{v}^{n + 1} - v^n}{\tau} =  {\Delta}_h \wh{u}^{n+\fr} -  \phi \sin{\ol{u}^{n+\fr}},\\[0.3cm]
		\dfrac{\wt{u}^{n+1} - u^n}{\tau} = \dfrac{\wt{v}^{n+1} + v^n}{2},\\[0.3cm]
		\dfrac{\wt{v}^{n+1} - v^n}{\tau} = {\Delta}_h \dfrac{\wt{u}^{n+1} + u^n}{2} -  \phi \sin {\wh{u}^{n+\fr} }.
		\end{cases}
		\end{align}
		\item We update $u^{n+1}, ~v^{n+1}$ by the following projection step \cite{hairer2006geometric}
		\begin{align}
		\begin{cases}
		u^{n+1} = \wt{u}^{n+1} + \lambda (- \Delta_h \wt{u}^{n+1} + \phi \sin{\wt{u}^{n+1}}),\\[0.3cm]
		v^{n+1} = \wt{v}^{n+1} + \lambda \wt{v}^{n+1},\\[0.3cm]
		\cH_h(u^{n+1},v^{n+1}) = \cH_h(u^0,v^0),
		\end{cases}
		\end{align}
		where $\lambda$ is a Lagrange multiplier and $\cH_h$ is the discrete energy given by
		\begin{align}
		\cH_h(u,v) = \frac{1}{2}\|v\|_h^2 - \frac{1}{2}(u,\Delta_h u)_h + (\phi, 1-\cos{u})_h.
		\end{align}	
	\end{enumerate}
\end{sch}

\begin{sch}[{\bf PEPM-R}]\label{energy-preserving-projection-PEPM-R}
	Given $u^{n-1}, ~u^n, ~v^n \in V_{h}^v, ~ \forall\; n \geq 1,$ we obtain $u^{n+1}, ~v^{n+1} \in V_{h}^v$ through the following two steps:
	\begin{enumerate}
		\item 
		We first compute $\wt{u}^{n+1}$ and $\wt{v}^{n+1}$ via the prediction-correction Crank-Nicolson scheme
		\begin{align}
		\begin{cases}
		\dfrac{\wh{u}^{n + 1} - u^n}{\tau} = \wh{v}^{n+\fr},\\[0.3cm]
		\dfrac{\wh{v}^{n + 1} - v^n}{\tau} =  \overline{\Delta}_h \wh{u}^{n+\fr} -  \phi \sin{\ol{u}^{n+\fr}},\\[0.3cm]
		\dfrac{\wt{u}^{n+1} - u^n}{\tau} = \dfrac{\wt{v}^{n+1} + v^n}{2},\\[0.3cm]
		\dfrac{\wt{v}^{n+1} - v^n}{\tau} = \overline{\Delta}_h \dfrac{\wt{u}^{n+1} + u^n}{2} -  \phi \sin {\wh{u}^{n+\fr} }.
		\end{cases}
		\end{align}
		\item We update $u^{n+1}, ~v^{n+1}$ by the following projection step
		\begin{align}
		\begin{cases}
		u^{n+1} = \wt{u}^{n+1} + \lambda (- \overline{\Delta}_h \wt{u}^{n+1} + \phi \sin{\wt{u}^{n+1}}),\\[0.3cm]
		v^{n+1} = \wt{v}^{n+1} + \lambda \wt{v}^{n+1},\\[0.3cm]
		\cH_h(u^{n+1},v^{n+1}) = \cH_h(u^0,v^0),
		\end{cases}
		\end{align}
		where
		\begin{align}
		\cH_h(u,v) = \frac{1}{2}\|v\|_T^2 - \frac{1}{2} \left<u, \overline{\Delta}_h u \right>_T + \left<\phi, 1-\cos{u} \right>_T.
		\end{align}	
	\end{enumerate}
\end{sch}

\begin{remark}
	Note that $\wt{u}^{n+1}$ and $\wt{v}^{n+1}$ in the above proposed schemes can be computed explicitly by DCT.	In the second step of the schemes, we can eliminate $u^{n+1}, ~v^{n+1}$ to derive a nonlinear algebraic equation for $\lambda$, which can be solved efficiently by the Newton iteration with 0 as the initial condition. 
\end{remark}

\section{Supplementary variable method (SVM)}\label{sect:SVM}
In this section, we propose a new perspective to develop energy-preserving algorithms for the SG equation. Firstly, we impose the energy conservation law as a constraint for the SG equation \eqref{model-eq}, i.e., 
\begin{align}\label{energy_constraint}
\cH[u(\bx,t),v(\bx,t)] \equiv \cH[u(\bx,0),v(\bx,0)], \quad \cH[u,v]  = \dfrac{1}{2} \int_{\Omega} \big( v^2 + | \nabla u|^2 + 2\phi( 1 - \cos{u}) \big) \d \bx.
\end{align}
Note that \eqref{model-eq} with the constraint \eqref{energy_constraint} constitutes an over-determined system. 
To solve them, we modify \eqref{model-eq} by a time-dependent supplementary variable $\beta(t)$ together with a user supplied function $g[u,v]$:
\begin{align}\label{SG-Lagrange-model}
\begin{cases}
u_t = v,\\
v_t = \Delta u - \phi \sin{u} + \beta(t)g[u, v],\\
\cH[u(\bx,t),v(\bx,t)] = \cH[u(\bx,0),v(\bx,0)],
\end{cases}
\end{align}
where $g[u, v]$ is a given function that may depend on $u$, $v$ and their derivatives. There is a grate deal of flexibility in determining how to relax the SG model with a supplementary variable. It is clearly an open problem for this approach. Here we choose $g[u,v] = \phi \sin{u}$ in this paper.

Applying the two cosine pseudo-spectral discretizations in space and the prediction-correction Crank-Nicolson scheme in time for the system \eqref{SG-Lagrange-model}, we obtain two new energy-preserving methods, which are named {\bf SVM-M} and {\bf SVM-R}, respectively.

\begin{sch}[{\bf SVM-M}]\label{fully-disc-SG-SVM-M}
	Given $u^{n-1}, ~u^n, ~v^n \in V_{h}^c, ~ \forall n \geq 1,$ we obtain $u^{n+1}, ~v^{n+1} \in V_{h}^c$ through the following two steps:
	\begin{enumerate}
		\item Prediction: predict $\wh{u}^{n + 1}$ and $\wh{v}^{n + 1}$ via an efficient and second-order schemes
		\begin{align}\label{fully-disc-SG-Pre}
		\begin{cases}
		\dfrac{\wh{u}^{n + 1} - u^n}{\tau} = \wh{v}^{n+\fr},\\
		\dfrac{\wh{v}^{n + 1} - v^n}{\tau} = \Delta_h \wh{u}^{n+\fr} -  \phi \sin{\ol{u}^{n+\fr}}.
		\end{cases}
		\end{align}
		\item Correction: 
		\begin{align}\label{fully-disc-SG-SVM-M-Correction}
		\begin{cases}
		\delta^+_t u^n = v^{n+\fr},\\
		\delta^+_tv^n = \Delta_h u^{n+\fr} -  \phi \sin { \wh{u}^{n+\fr} }  + \beta^{n+\fr} g [\wh{u}^{n+\fr}, \wh{v}^{n+\fr}],\\
		\cH_h[u^{n+1},v^{n+1}] = \cH_h[u^0, v^0],
		\end{cases}
		\end{align}
		where \begin{align}
		\cH_h[u,v] = \frac{1}{2}\|v\|_h^2 - \frac{1}{2}(u,\Delta_h u)_h + (\phi, 1-\cos{u})_h.
		\end{align}	
	\end{enumerate}
\end{sch}

\begin{sch}[{\bf SVM-R}]\label{fully-disc-SG-newSAV-model-N-N}
	Given $u^{n-1}, ~u^n, ~v^n \in V_{h}^v, ~ \forall n \geq 1,$ we obtain $u^{n+1}, ~v^{n+1} \in V_{h}^v$ by the following two steps:
	\begin{enumerate}
		\item Prediction: predict $\wh{u}^{n + 1}$ and $\wh{v}^{n + 1}$ via an efficient and second-order schemes
		\begin{align}\label{fully-disc-SG-newSAV-model-N-N-Pre}
		\begin{cases}
		\dfrac{\wh{u}^{n + 1} - u^n}{\tau} = \wh{v}^{n+\fr},\\
		\dfrac{\wh{v}^{n + 1} - v^n}{\tau} = \overline{\Delta}_h \wh{u}^{n+\fr} -  \phi \sin{\ol{u}^{n+\fr}}.
		\end{cases}
		\end{align}
		\item Correction: 
		\begin{align}\label{fully-disc-SG-newSAV-model-N-N-Cor}
		\begin{cases}
		\delta^+_t u^n = v^{n+\fr},\\
		\delta^+_tv^n = \overline{\Delta}_h u^{n+\fr} -  \phi \sin { \wh{u}^{n+\fr} }  + \beta^{n+\fr} g [\wh{u}^{n+\fr}, \wh{v}^{n+\fr}],\\
		\cH_h[u^{n+1},v^{n+1}] = \cH_h[u^0, v^0],
		\end{cases}
		\end{align}
		where
		\begin{align}
		\cH_h[u,v] = \frac{1}{2}\|v\|_T^2 - \frac{1}{2} \left<u, \overline{\Delta}_h u \right>_T + \left<\phi,  1-\cos{u} \right>_T.
		\end{align}	
	\end{enumerate}
\end{sch}

In the following, we show how to solve {\bf SVM-R} efficiently. Note that {\bf SVM-M} can be solved similarly. According to \eqref{fully-disc-SG-newSAV-model-N-N-Pre}, we have 
\begin{align}
\wh{u}^{n+\fr} = (\bI -\frac{\tau^2}{4} \overline{\Delta}_h)^{-1} (u^n+\frac{\tau}{2} v^n - \frac{\tau^2}{4} \phi\sin{\ol{u}^{n+\fr}}), \quad \wh{v}^{n+\fr} = \frac{2}{\tau} (\wh{u}^{n+\fr} - u^n).
\end{align}	
Letting 
\begin{align*}
& \wt{u}^{n+1} = (\bI - \frac{\tau^2}{4}\overline{\Delta}_h)^{-1}\left( (\bI + \frac{\tau^2}{4}\overline{\Delta}_h)u^n + \tau v^n - \frac{\tau^2}{2}   \phi \sin{\wh{u}^{n+\fr}}  \right), \quad \wt{v}^{n+1} = \frac{2}{\tau}(\wt{u}^{n+1} - u^n ) - v^n, \\
& \omega^{n} =  \dfrac{\tau^2}{2}(\bI - \frac{\tau^2}{4}\overline{\Delta}_h)^{-1}g[\wh{u}^{n+
	\fr}, \wh{v}^{n+\fr}], \quad \gamma^n = \frac{2}{\tau} \omega^{n},
\end{align*}
we can deduce from \eqref{fully-disc-SG-newSAV-model-N-N-Cor} that
\begin{numcases}{}
u^{n+1} = \wt{u}^{n+1} + \beta^{n+\fr} \omega^n, \label{lineq-solver-eq2} \\
v^{n+1} = \wt{v}^{n+1} + \beta^{n+\fr} \gamma^n. \label{lineq-solver-eq3}
\end{numcases}
Then, substituting  \eqref{lineq-solver-eq2}-\eqref{lineq-solver-eq3} into the last equation of \eqref{fully-disc-SG-newSAV-model-N-N-Cor} leads to
\begin{align}\label{nonlinear-eta_t-eq}
\cH_h[\wt{u}^{n+1} + \beta^{n+
	\fr}\omega^n, \wt{v}^{n+1} + \beta^{n+\fr} \gamma^n] = \cH_h[u^0, v^0],
\end{align}
which is a scalar nonlinear equation for $\beta^{n+\fr}$. In general, it can have multiple solutions, but one of them must approximate to zero as $\tau \rightarrow 0$. Therefore, we solve for this solution by using an iterative method such as Newton iteration with $0$ as the initial value, it generally converges to a solution close to $0$ when $\tau$ is not too large. After obtaining $\beta^{n+\fr}$, we update $u^{n+1}$ and $v^{n+1}$ using \eqref{lineq-solver-eq2} and \eqref{lineq-solver-eq3}, respectively. 

\begin{remark}
	The SVM can be regarded as a perturbation or projection to the PDE model. With proper choice of discretization and the way that supplementary variable is added to the original  model, it is similar to the first approach in the process of fast solving SVM. 
	The supplementary function $g[u, v]$  in SVM is more flexible, and the numerical experiment also verifies the fact.  
	The idea of SVM has been applied  to deal with dissipative system \cite{Hong2020&AML,Cheng2019,Gong2020&SVM}. Philosophically, SVM originates from a quite different perspective to study  the development of structure-preserving schemes.
\end{remark}

\section{Numerical Results}\label{sec:NR}
In this section, we perform several numerical experiments to confirm the convergence rate of the proposed schemes and to study the energy conservative behavior for the sine-Gordon model subject to the Neumann boundary conditions at mid-point and regular grid. In numerical implement, replacing $\ol{u}^{n+\fr}$ with $u^n$ in the prediction step leads to a two-level scheme, which is used to compute the initial value. Unless otherwise stated, the default value of the iteration tolerance is set as $\text{Tol} = 1.0e-14$.

\subsection{Accuracy test and numerical comparisons on accuracy of solution}
\begin{example}[Mesh refinement test]\label{eg:MRT}
	We consider the one-dimensional sine-Gordon equation 
	\begin{align*}
	u_{tt} - u_{xx} + \sin{u} = 0,
	\end{align*}
	which has a theoretical solution \cite{bratsos2008numerical}
	\begin{align*}
	u(x, t) = 4 \arctan \left( c^{-1} \sin(c \kappa t) \sech(\kappa x) \right),
	\end{align*}
	where $\kappa = 1/\sqrt{1 + c^2}$. 
	The initial conditions are given by
	\begin{align*}
	u(x, 0) = 0,\quad
	u_t(x, 0) = 4 \kappa \sech(\kappa x).
	\end{align*}
	It is known as the breather solution of the sine-Gordon equation, which represents a pulse type soliton. The computational domain is set as $\Omega = [-20, 20]$ and the parameter $c$ is the velocity and $c = 0.5$. For the spatial test, we choose the time step as $\tau = 1.0e-4$ to prevent the errors in time discretization from contaminating our results. With grid sizes from $N_x = N =16, 32, 64, 128$ to  $256$, the errors in the $L^2$ and $L^{\infty}$ are calculated up to time $t = 1$. These results are summarized in Figure \ref{fig:test-space}, where we observe the spectral accuracy in space for the four schemes. 
	\begin{figure}[H]
		\centering
		\subfigure[The $L^2$ error plot for $u$ in space.]{
			\includegraphics[width=0.4\textwidth,height=0.25\textwidth]{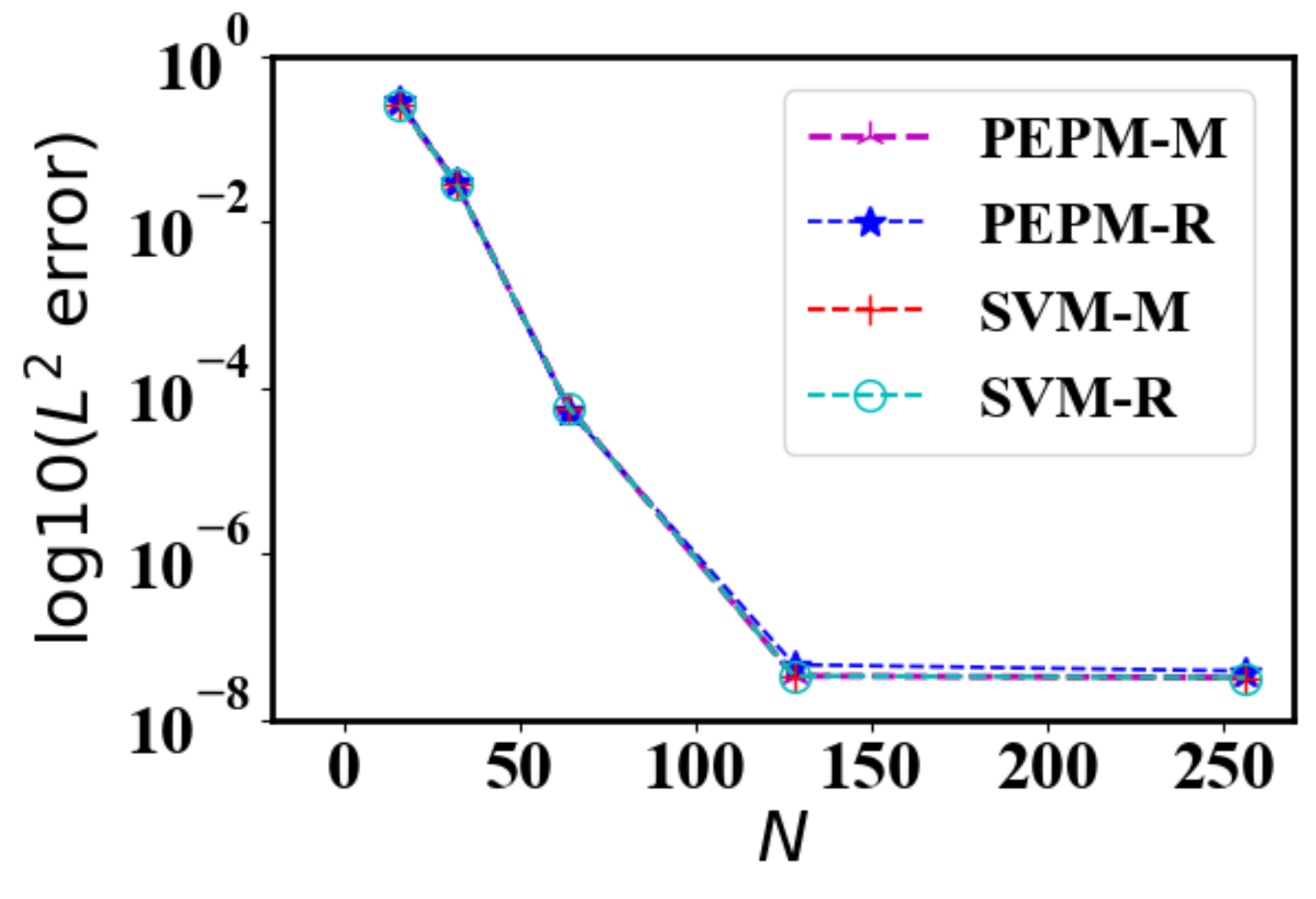}
		}\quad
		\subfigure[The $L^{\infty}$ error plot for $u$ in space.]{
			\includegraphics[width=0.4\textwidth,height=0.25\textwidth]{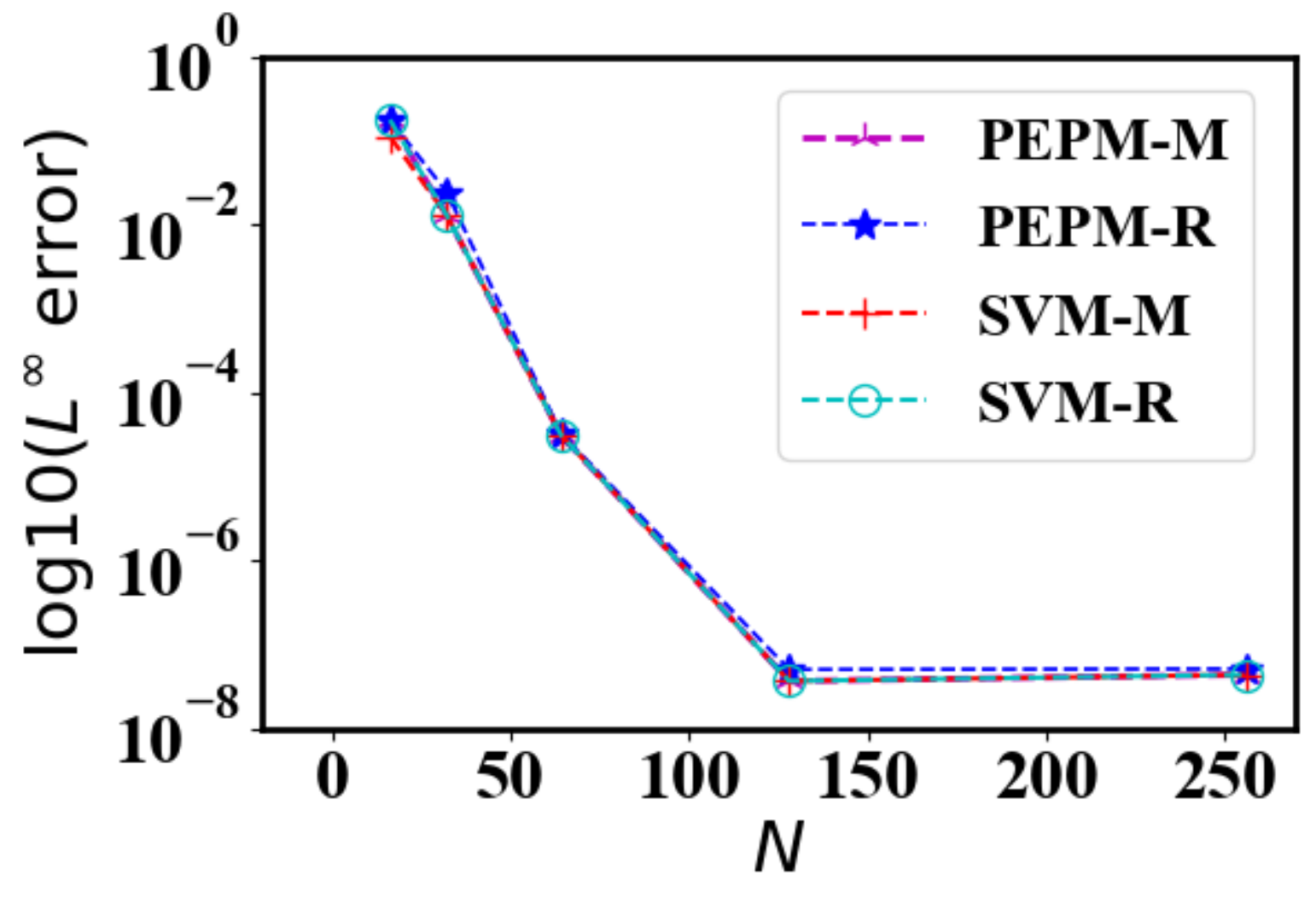}
		}
		\caption{{\bf Example \ref{eg:MRT}}: Mesh refinement test for space accuracy at a small time step $\tau = 1.0e-4$. A spectral accuracy is observed.} \label{fig:test-space}
	\end{figure}
	For time accuracy test, we fix $N = 256$. The discrete  $L^2$ and $L^{\infty}$ errors in time for four schemes at the final time $T = 1$ are summarized in Figure \ref{fig:test-order1}. It is observed that the four schemes achieve the expected second order convergence in time. Moreover,  Figure \ref{fig:test-order2} shows the second-order  and third-order accuracy for the supplementary variable $\beta$ and the Lagrange multiplier $\lambda$ are reached, respectively. Thus, the numerical performances from Example \ref{eg:MRT} validate the correctness of our proposed schemes.
	\begin{figure}[H]
		\centering
		\subfigure[The $L^2$ error plot for $u$ in time.]{
			\includegraphics[width=0.45\textwidth,height=0.225\textwidth]{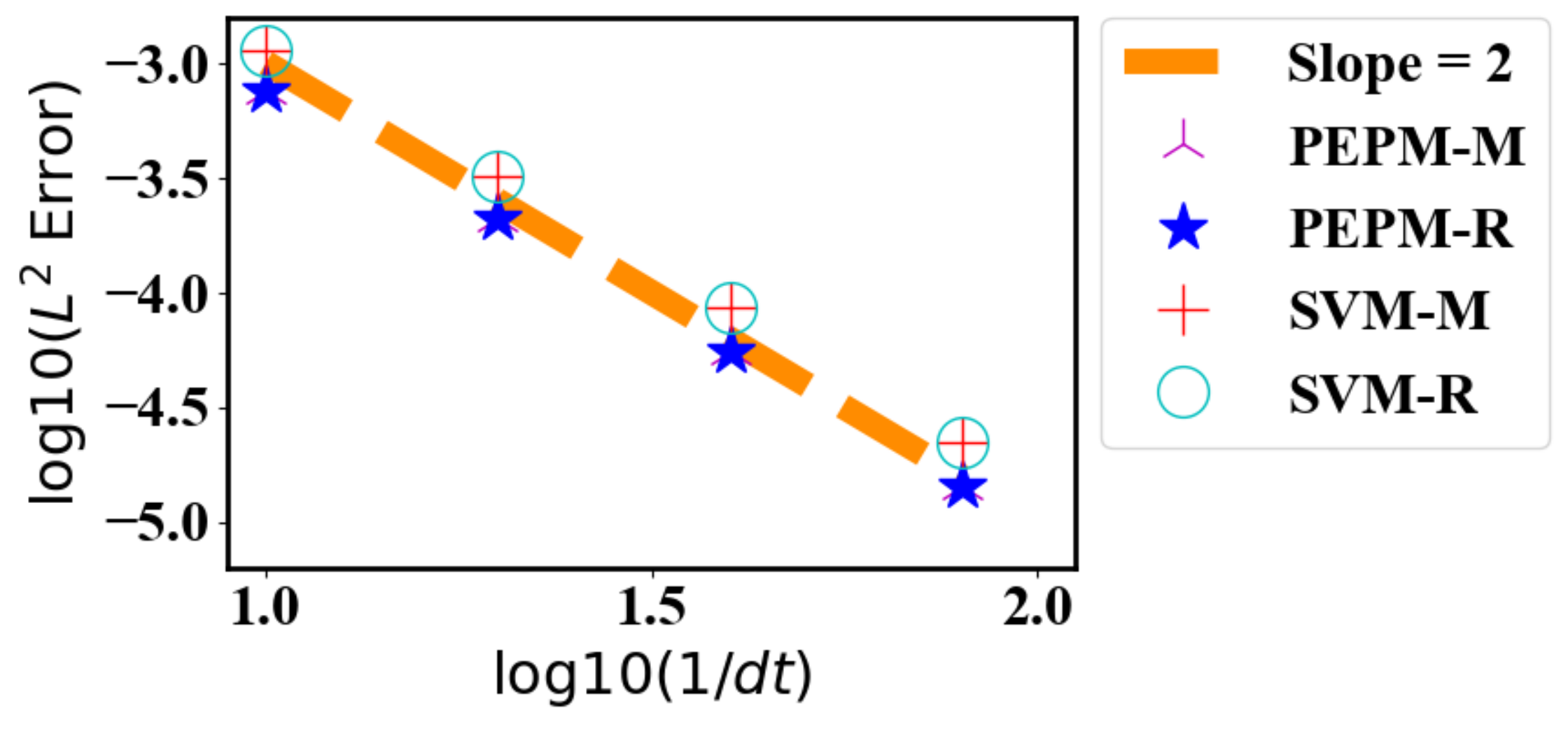}
		}
		\subfigure[The $L^{\infty}$ error plot for $u$ in  time.]{
			\includegraphics[width=0.45\textwidth,height=0.225\textwidth]{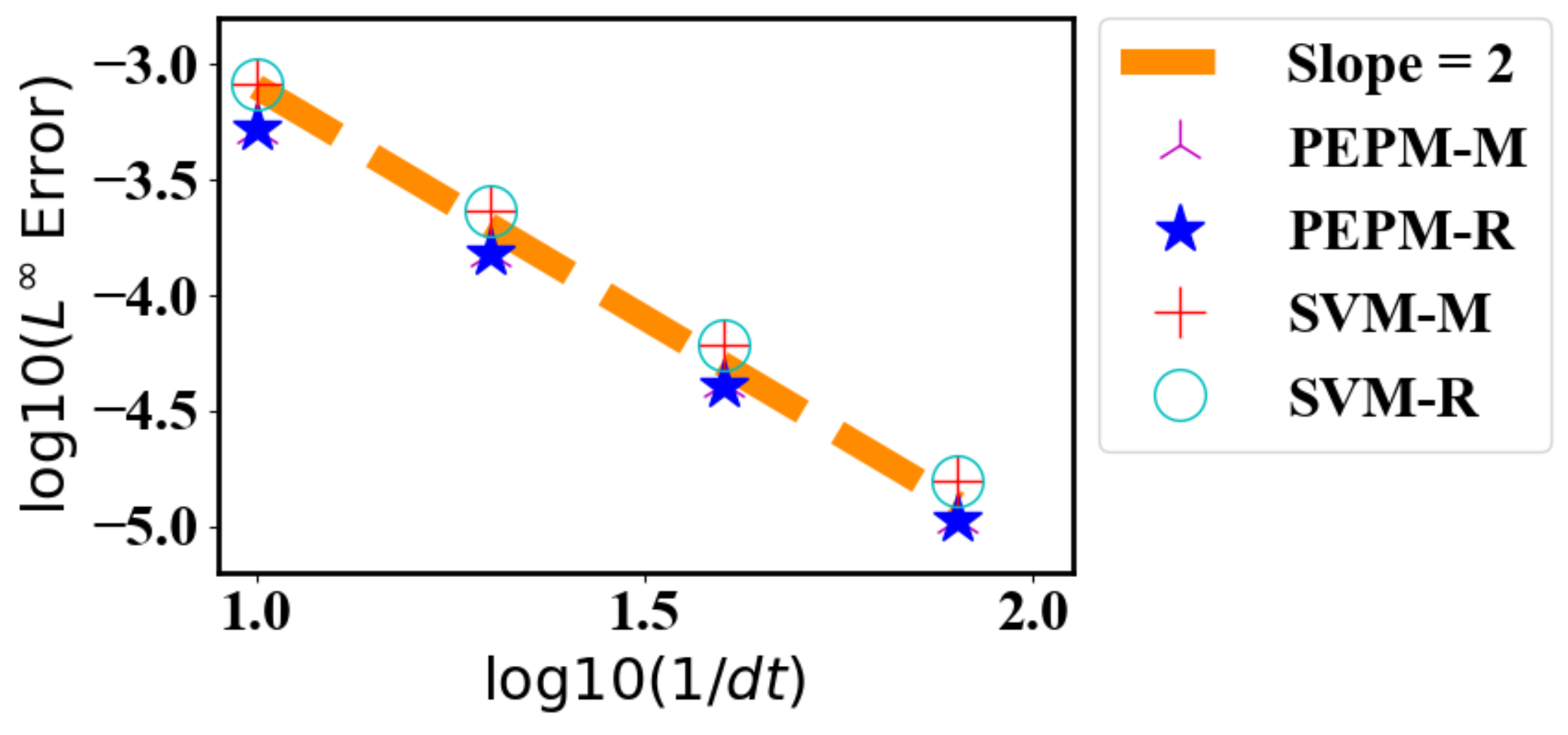}
		}
		\caption{{\bf Example \ref{eg:MRT}}: Mesh refinement test in time for four different schemes. Here, we fix spatial meshes at $N = 256$. Desired order of numerical solution accuracy (2th order) is achieved.} \label{fig:test-order1}
		
	\end{figure}
	\begin{figure}[H]
		\centering
		\subfigure{
			\includegraphics[width=0.55\textwidth,height=0.25\textwidth]{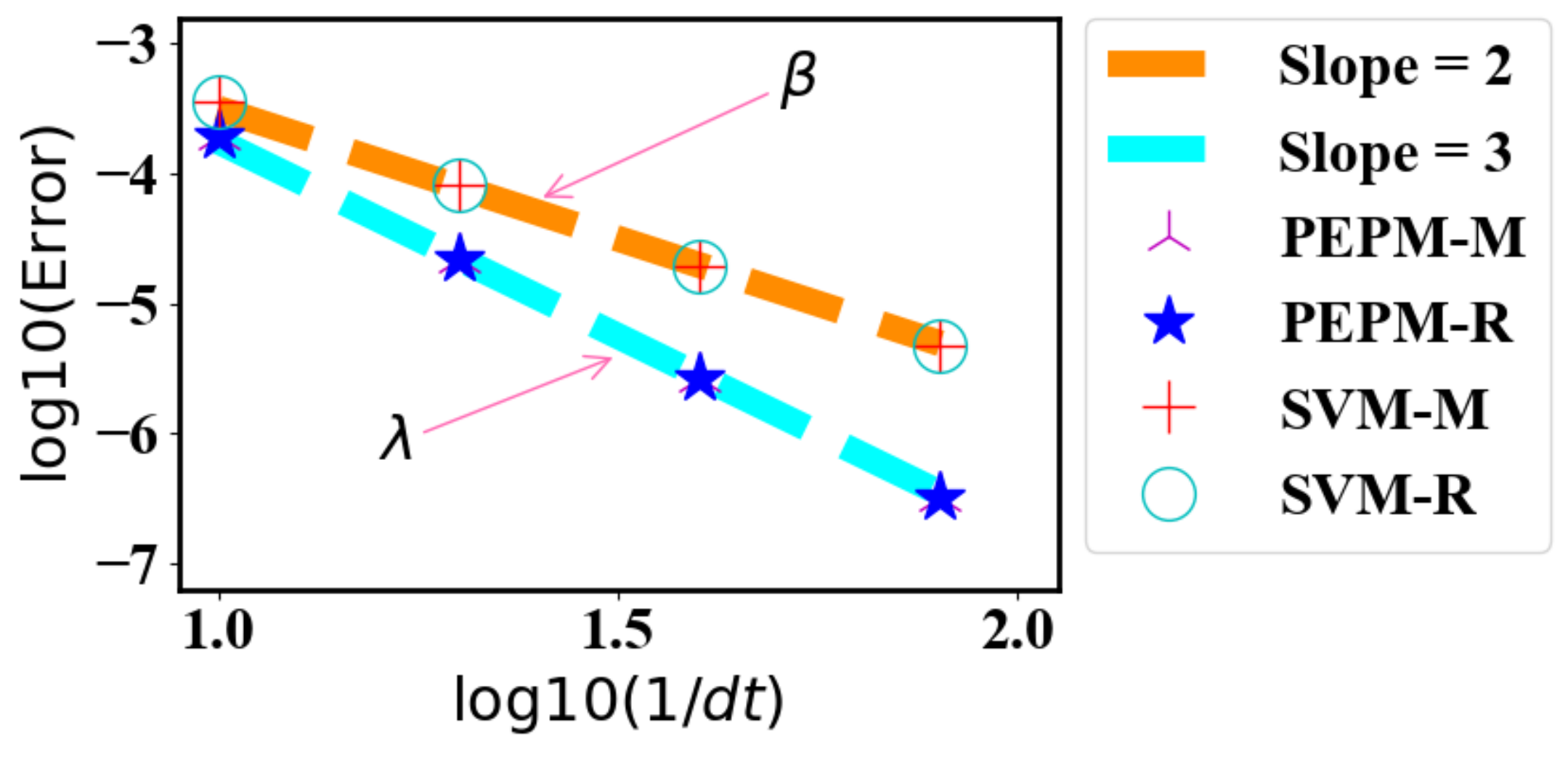}
		}
		\caption{{\bf Example \ref{eg:MRT}}:  The error plots of supplementary variable and Lagrange multiplier. The slopes of  $\beta$ and $\lambda$ error curves are close to $2$ and $3$, respectively.} \label{fig:test-order2}
	\end{figure}
	Furthermore, we make some numerical comparisons on accuracy of solution between our proposed schemes and the energy-preserving schemes based on the SAV approach \cite{cai2019structure} $(${\bf SAV-M}, {\bf SAV-R} for short$)$, and the  AVF energy-preserving methods \cite{celledoni2012preserving} $(${\bf AVF-M}, {\bf AVF-R} for short$)$, where the cosine pseudo-spectral methods are used for spatial discretization on mid-point and regular grid in all numerical methods. The corresponding discrete $L^2$ errors and the maximum errors in solution are displayed in Table \ref{tab:Numer-Compari-solution}.
	These numerical results demonstrate that  our proposed four numerical schemes can provide more accurate solution than  the others. It is interesting that {\bf PEPM-M} exhibits slightly better than {\bf SVM-M} on mid-point grid.
	
	\begin{table}[H]
		{\caption{The errors in solution by various methods at $t = 10$ with $\tau = 1.0e-2$ and $N = 128$.}\label{tab:Numer-Compari-solution}}
		\begin{center}
			\begin{tabular}{c c c c c }	\toprule[1.0pt]
				Method  &{\bf SAV-M} & {\bf AVF-M}    &{\bf PEPM-M}   &{\bf SVM-M}  \\ \hline\specialrule{0em}{1pt}{1pt}
				$L^2$ error & 1.50e-01 &4.12e-05  &9.20e-06 &2.35e-05  \\ [0.1cm]
				$L^{\infty}$ error & 7.17e-02 &2.14e-05   &6.70e-06 &1.20e-05\\ \midrule
				Method  &{\bf SAV-R} & {\bf AVF-R}  &{\bf PEPM-R} &{\bf SVM-R}        \\ \hline\specialrule{0em}{1pt}{1pt}
				$L^2$ error & 1.50e-01 &4.13e-05 &3.35e-05 &2.35e-05 \\ [0.1cm]
				$L^{\infty}$ error & 7.19e-02 &2.08e-05&1.57e-05 &1.15e-05 \\ 
				\bottomrule[1.0pt]
			\end{tabular}
		\end{center}
	\end{table}
	
\end{example}
\subsection{Line solitons}
\begin{example}[Perturbation of a line soliton]\label{eg:PLS}
	In this example, we consider  the case of $\phi(x,y) = 1$ and the initial conditions
	\begin{align*}
	\begin{cases}
	u(\bx, 0) = 4 \tan^{-1}\left[ \exp(x + 1 - 2\sech (y + 7) - 2\sech( y - 7) ) \right],\\
	v(\bx, 0) = 0.
	\end{cases}
	\end{align*}
	In this case, we use the codes developed from the four schemes to simulate perturbation of a line soliton  in a domain $\Omega = [-7, 7]^2$ with $128\times 128$ meshes. We conduct this simulation with $\tau = 0.01$ and display  the isolines of $\sin(u/2)$ up to $T = 11$.  In Figure \ref{fig:Perturbation-evolution}, we observe that two symmetric dents moving toward each other, collapsing at $t = 7$ and continue to move away after the collision. These profiles look qualitatively similar to the reported results \cite{christiansen1981numerical}. Figure \ref{fig:Perturbation-energy-error} (a)  illustrates that our proposed schemes preserve the original energy conservation very accurately. Figure \ref{fig:Perturbation-energy-error} (b) depicts the supplementary variable $\beta $ and the Lagrange multiplier $\lambda $ oscillate around zero. 
	\begin{figure}[H]
		\centering
		\subfigure[$t=0$]{
			\includegraphics[width=0.22\textwidth,height=0.22\textwidth]{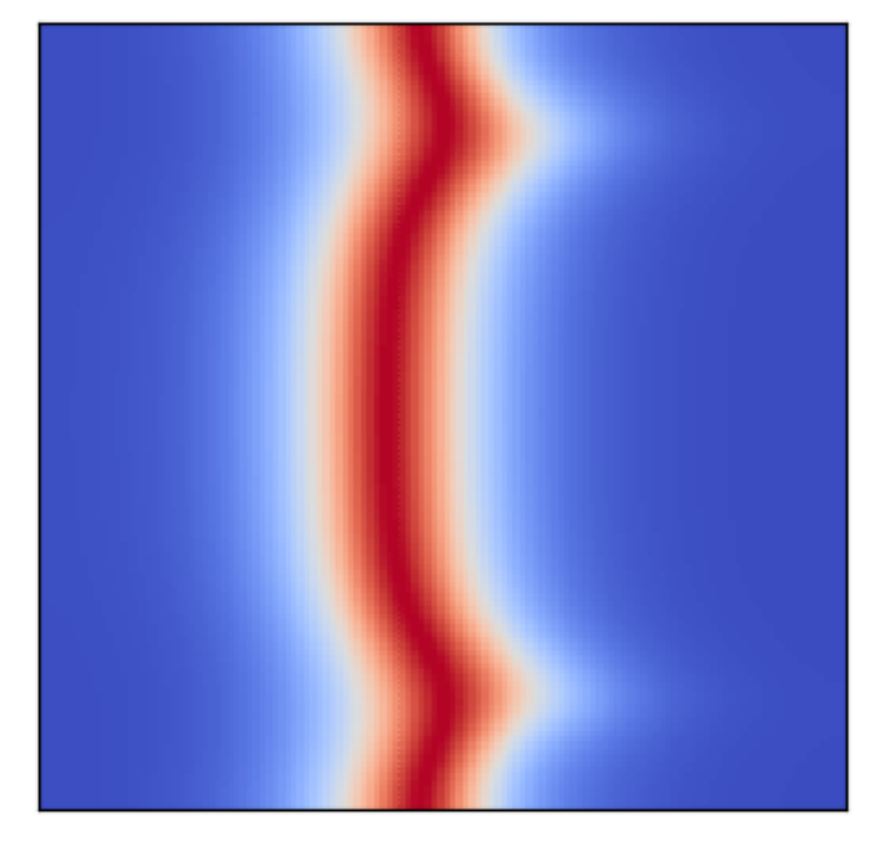}
		}
		\subfigure[$t=2$]{
			\includegraphics[width=0.22\textwidth,height=0.22\textwidth]{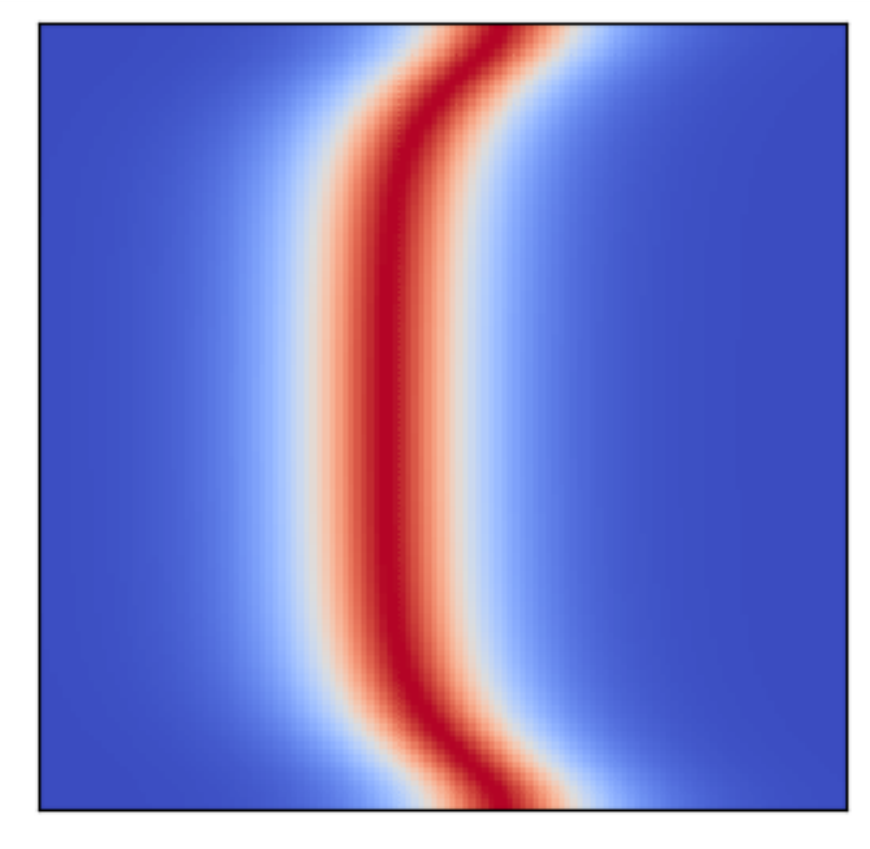}
		}
		\subfigure[$t=5$]{
			\includegraphics[width=0.22\textwidth,height=0.22\textwidth]{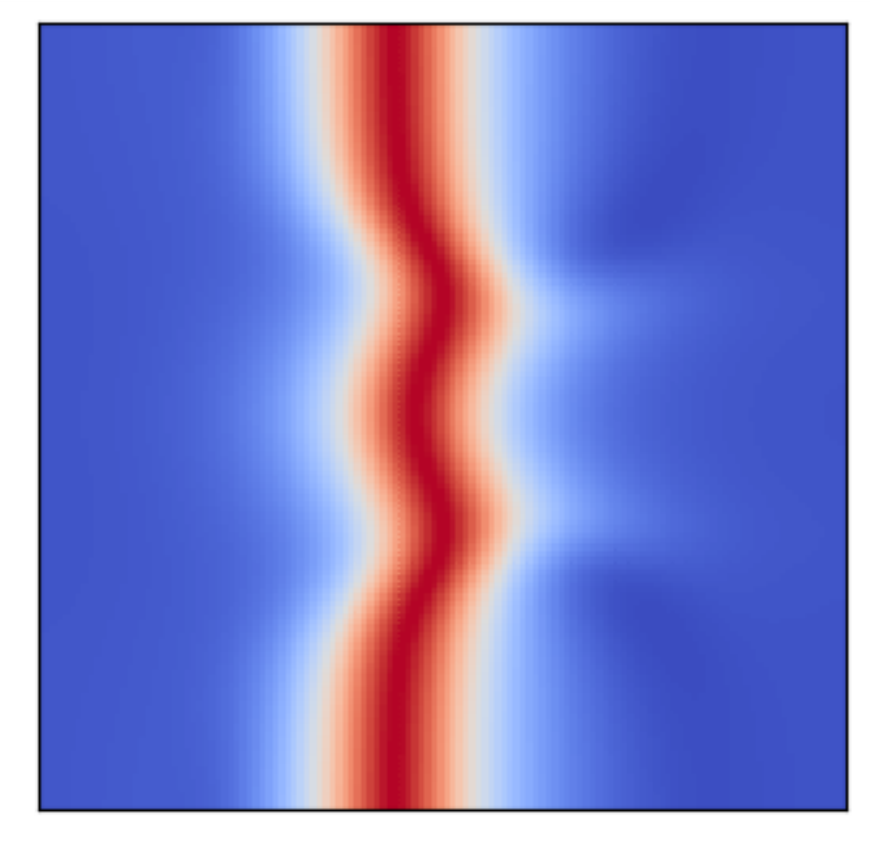}
		}\\~~
		\subfigure[$t=7$]{
			\includegraphics[width=0.22\textwidth,height=0.22\textwidth]{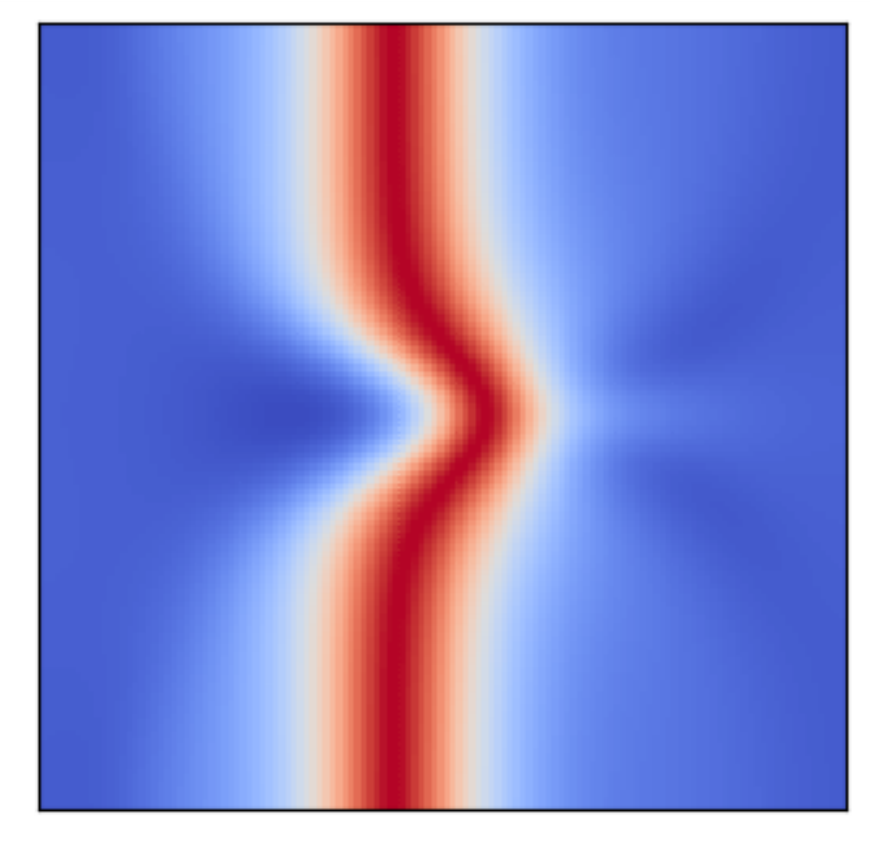}
		}
		\subfigure[$t=9$]{
			\includegraphics[width=0.22\textwidth,height=0.22\textwidth]{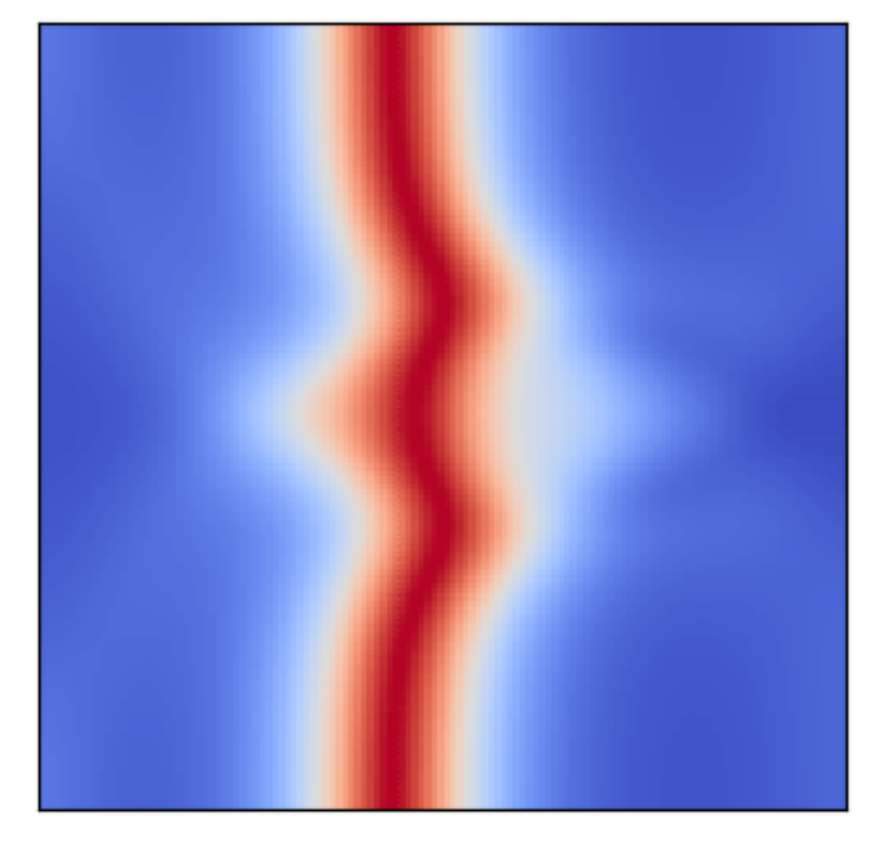}
		}
		\subfigure[$t=11$]{
			\includegraphics[width=0.22\textwidth,height=0.22\textwidth]{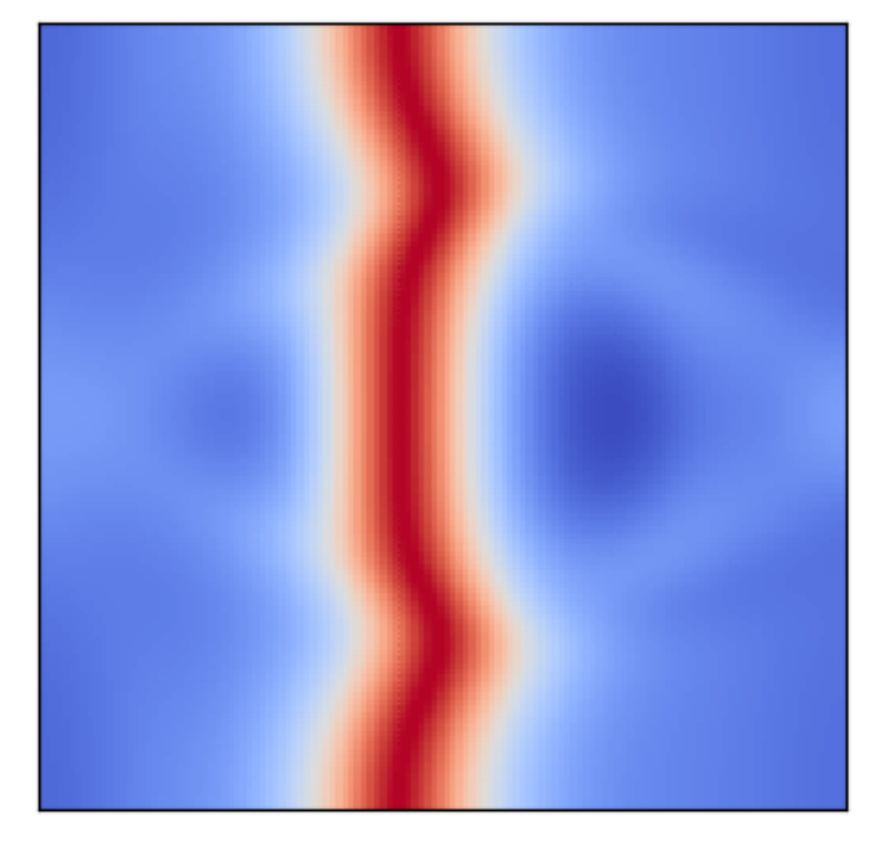}
		}
		\caption{{\bf Example \ref{eg:PLS}}: The isolines of numerical solutions of $\sin(u/2)$. Snapshots are taken at $t = 0, 2, 5, 7, 9, 11$, respectively.   \label{fig:Perturbation-evolution}}
	\end{figure}
	\begin{figure}[H]
		\centering
		\subfigure[]{
			\includegraphics[width=0.38\textwidth,height=0.28\textwidth]{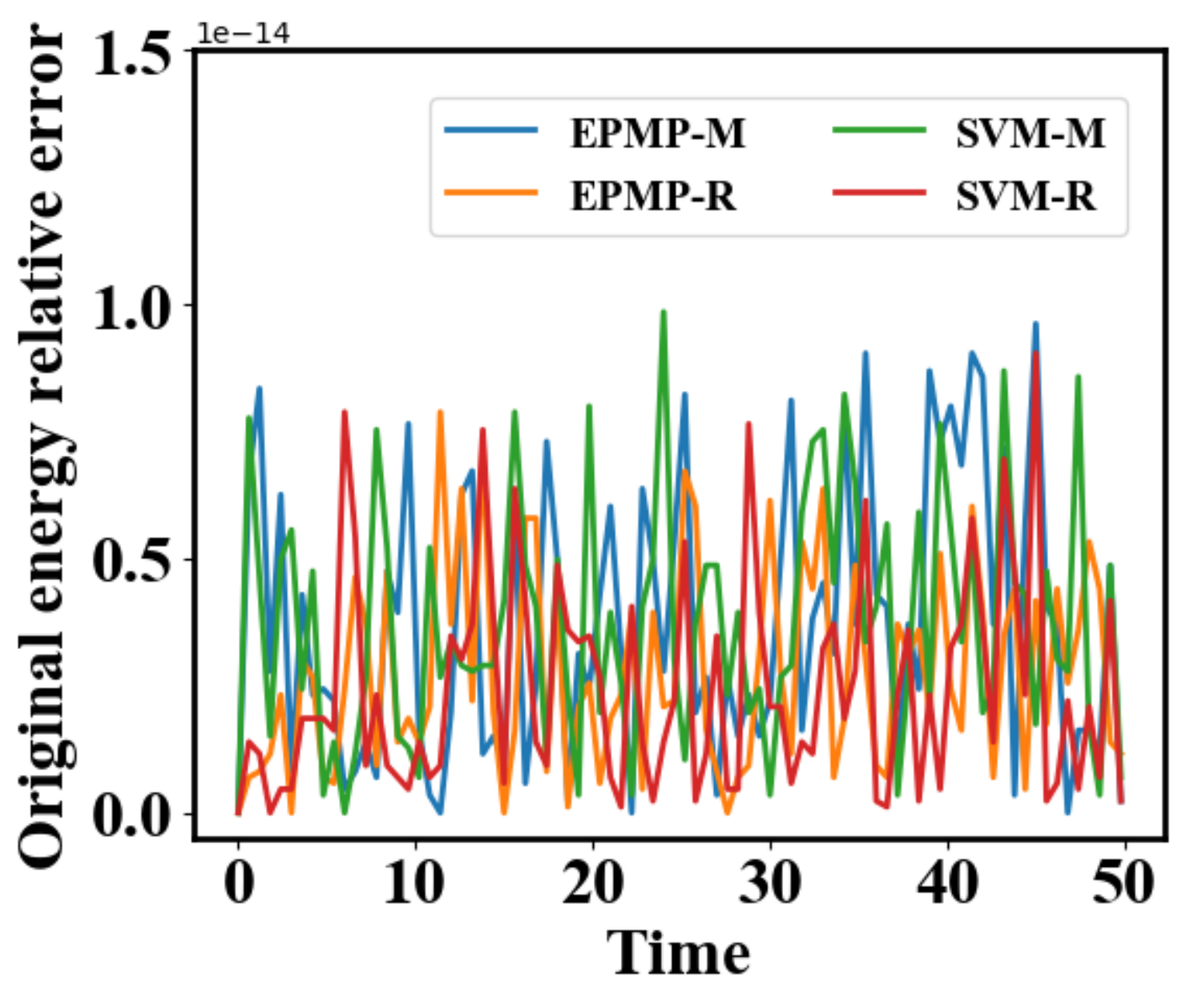}
		}
		\subfigure[]{
			\includegraphics[width=0.38\textwidth,height=0.28\textwidth]{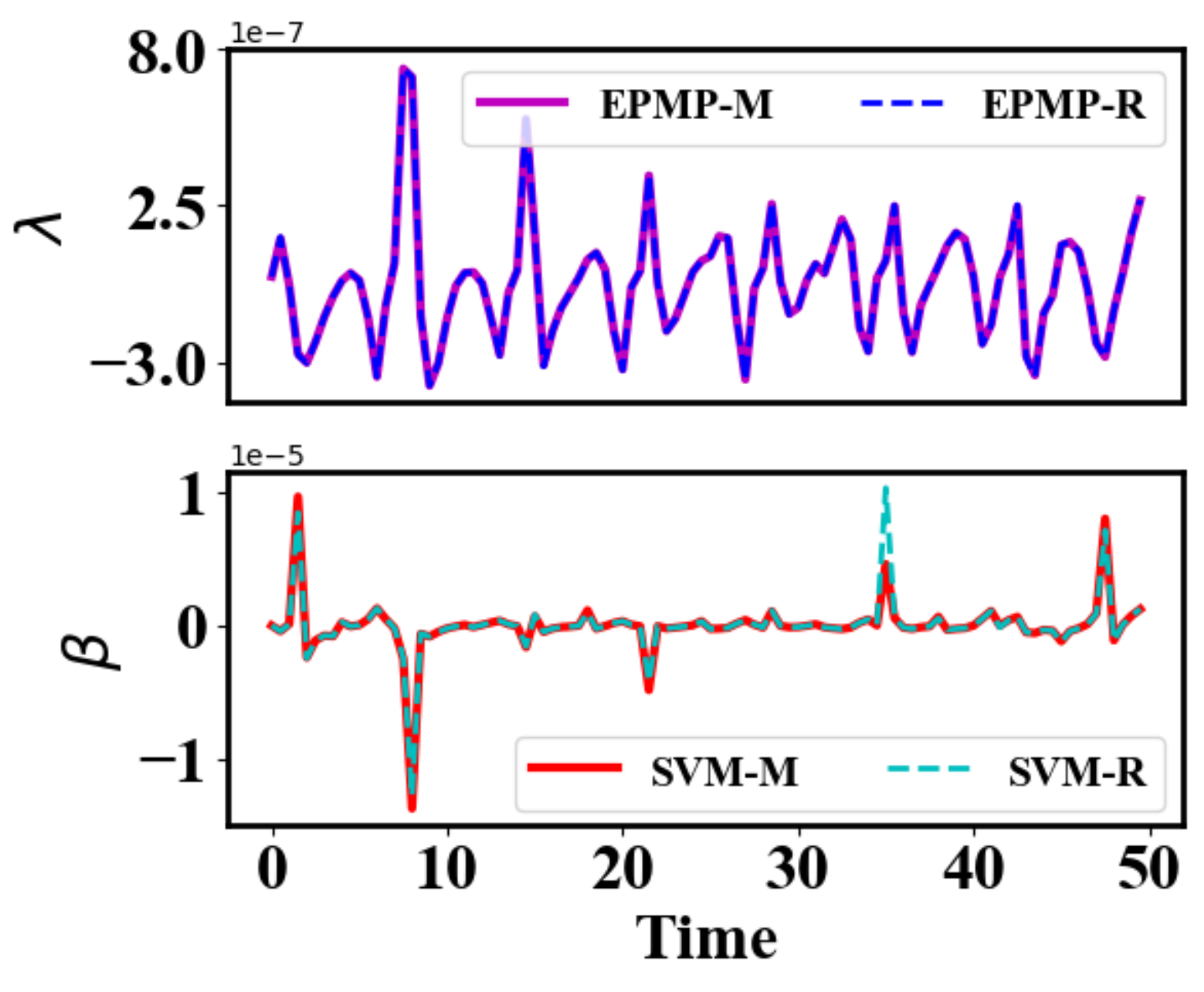}
		}
		\caption{{\bf Example \ref{eg:PLS}}: (a) Evolution of the  energy errors using the four schemes  with $\tau = 0.01$ and $N_x = N_y =128$. The curves of the energy errors show the four schemes warrant the energy conservation law very accurately with relatively large time step. (b) Evolution of the supplementary variable $\beta$ and the Lagrange multiplier $\lambda$.
			\label{fig:Perturbation-energy-error}}
	\end{figure}
	To further compare the advantages of our proposed schemes with   {\bf SAV-M}, {\bf SAV-R}, {\bf AVF-M} and {\bf AVF-R}, we summarize their  computational costs and energy errors in Figure \ref{fig:cpu-time} and Figure \ref{fig:different-energy-error}.
	In these simulations, we set $Tol = 1.0e-12$. 
	The datas from Figure \ref{fig:cpu-time}  support these observations that
	{\bf SVM-M/-R} and {\bf PEPM-M/-R} for this test are less efficient than 
	{\bf SAV-M/-R}, while much more efficient than the schemes {\bf AVF-M/-R}.
	The  price we pay using our proposed schemes in terms of CPU computing efficiency is that we have to solve a scalar  nonlinear equation whose cost is negligible compared to the prediction-correction step, while {\bf AVF-M/-R} solves a nonlinear system.
	Figure \ref{fig:different-energy-error} shows that our proposed schemes admit the original energy conservation law very well, but  fails to {\bf SAV-M/-R}.
	
	\begin{figure}[H]
		\centering
		\subfigure[On mid-point gird]{
			\includegraphics[width=0.40\textwidth,height=0.30\textwidth]{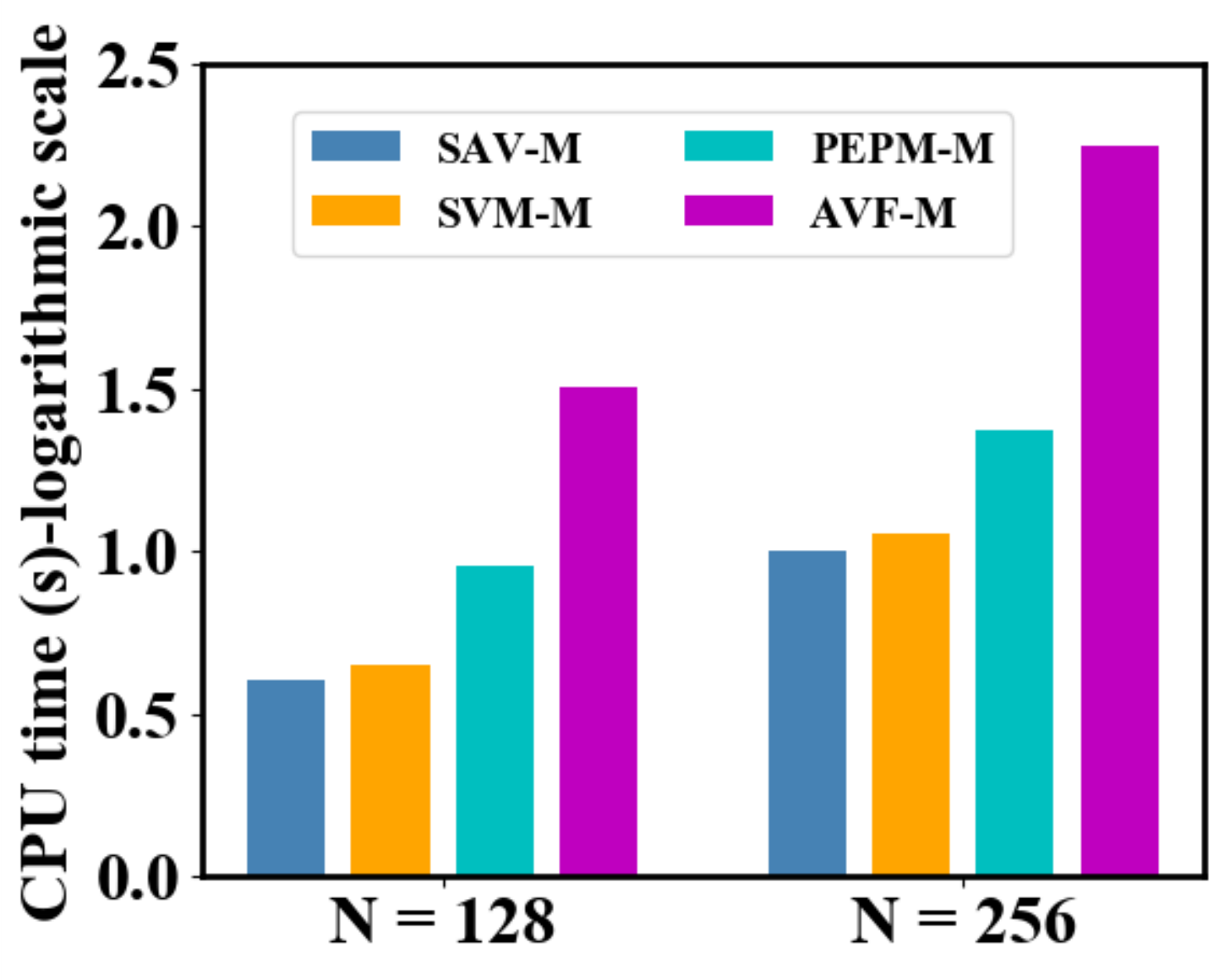}
		}
		\subfigure[On regular grid]{
			\includegraphics[width=0.40\textwidth,height=0.30\textwidth]{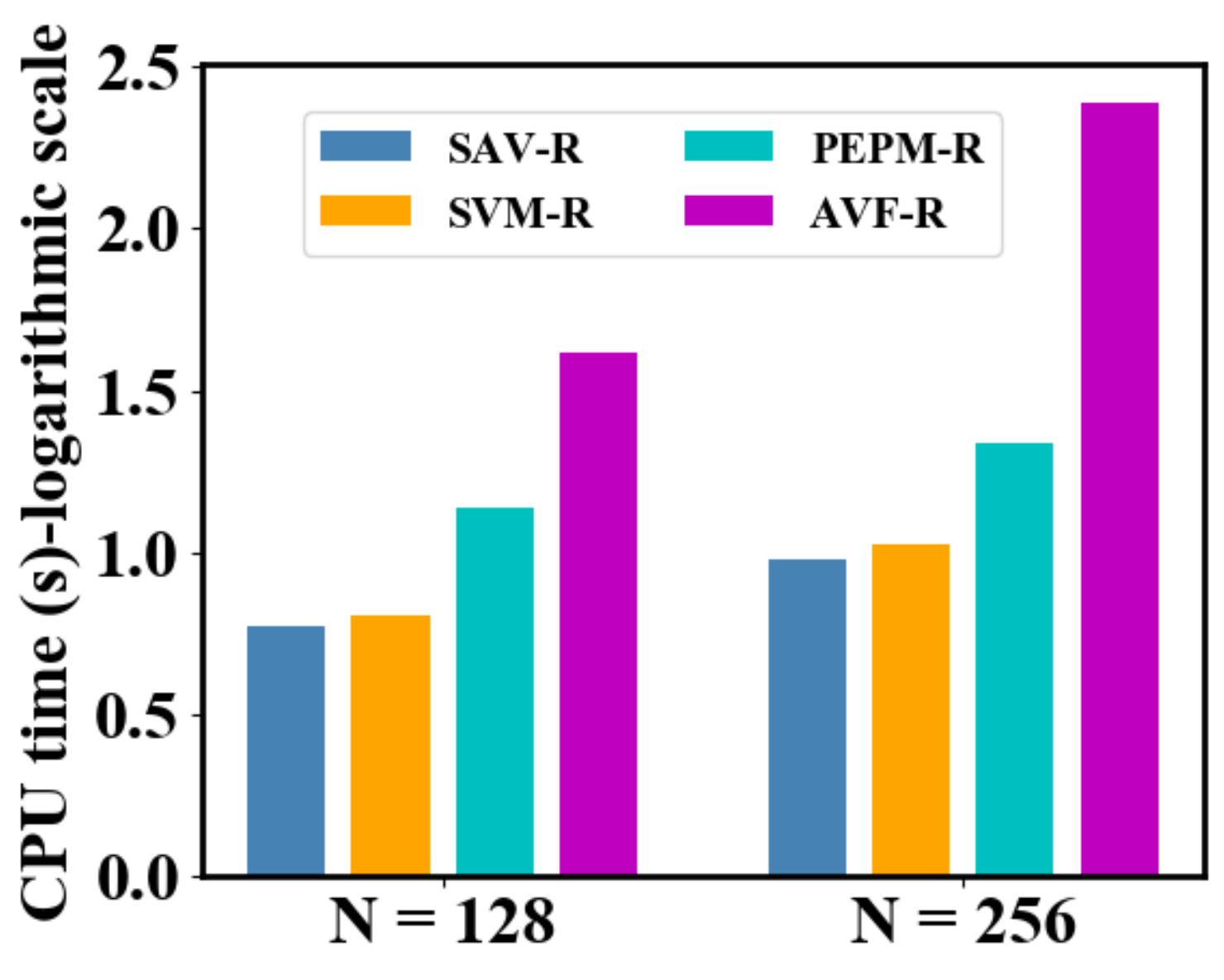}
		}
		\caption{{\bf Example \ref{eg:PLS}}: Comparison of CPU times in logarithmic scale  using different numerical methods with various spatial mesh sizes till $t = 1$, where the time step is set as $\tau = 0.01$. The bar charts show that our proposed schemes perform more slowly than  {\bf SAV-M/-R} in computational efficiencies, while present more superior than {\bf AVF-M/-R}.
			\label{fig:cpu-time}}
	\end{figure}
	
	\begin{figure}[H]
		\centering
		\subfigure[Simulation on mid-point gird]{
			\includegraphics[width=0.38\textwidth,height=0.28\textwidth]{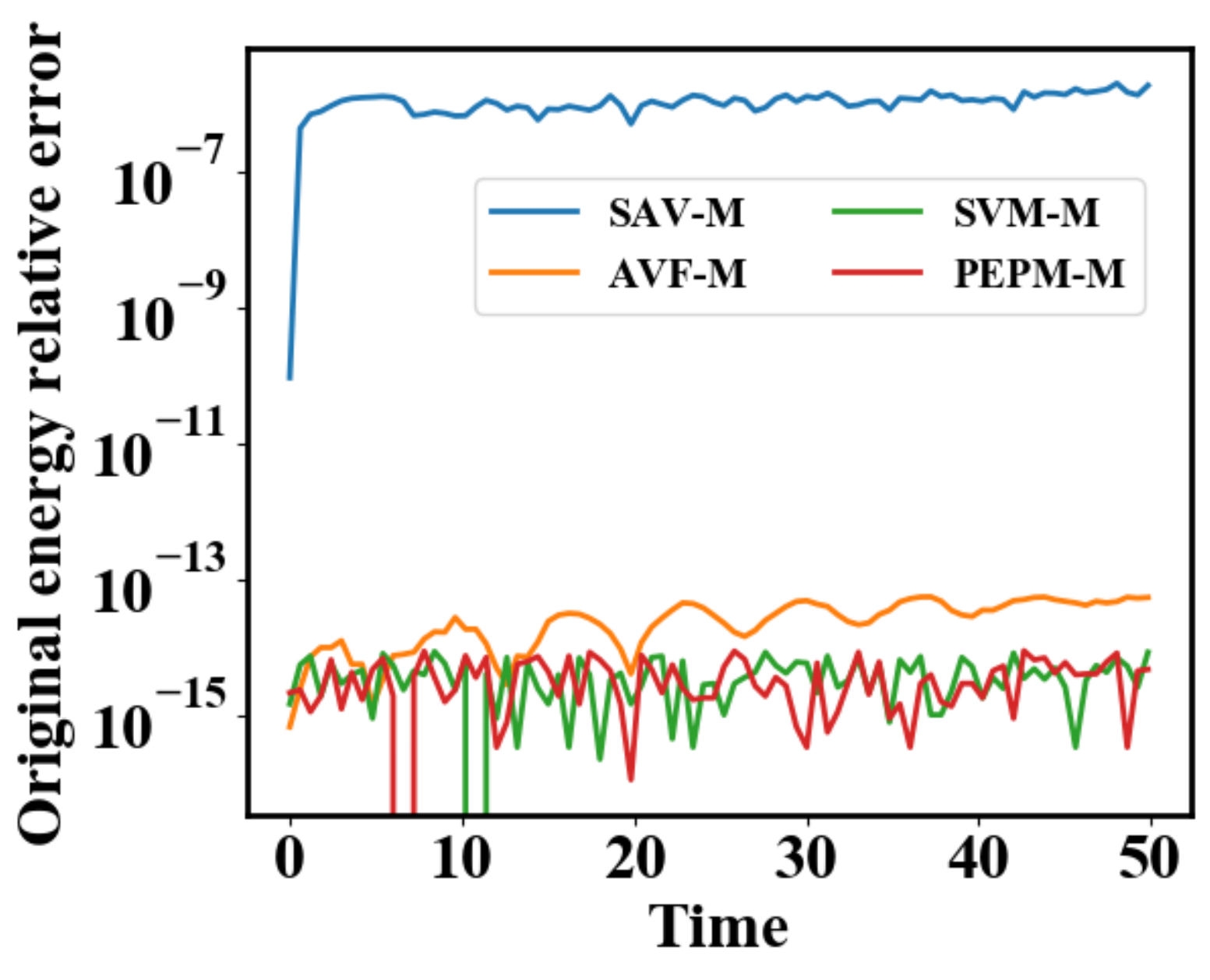}
		}
		\subfigure[Simulation on regular grid]{
			\includegraphics[width=0.38\textwidth,height=0.28\textwidth]{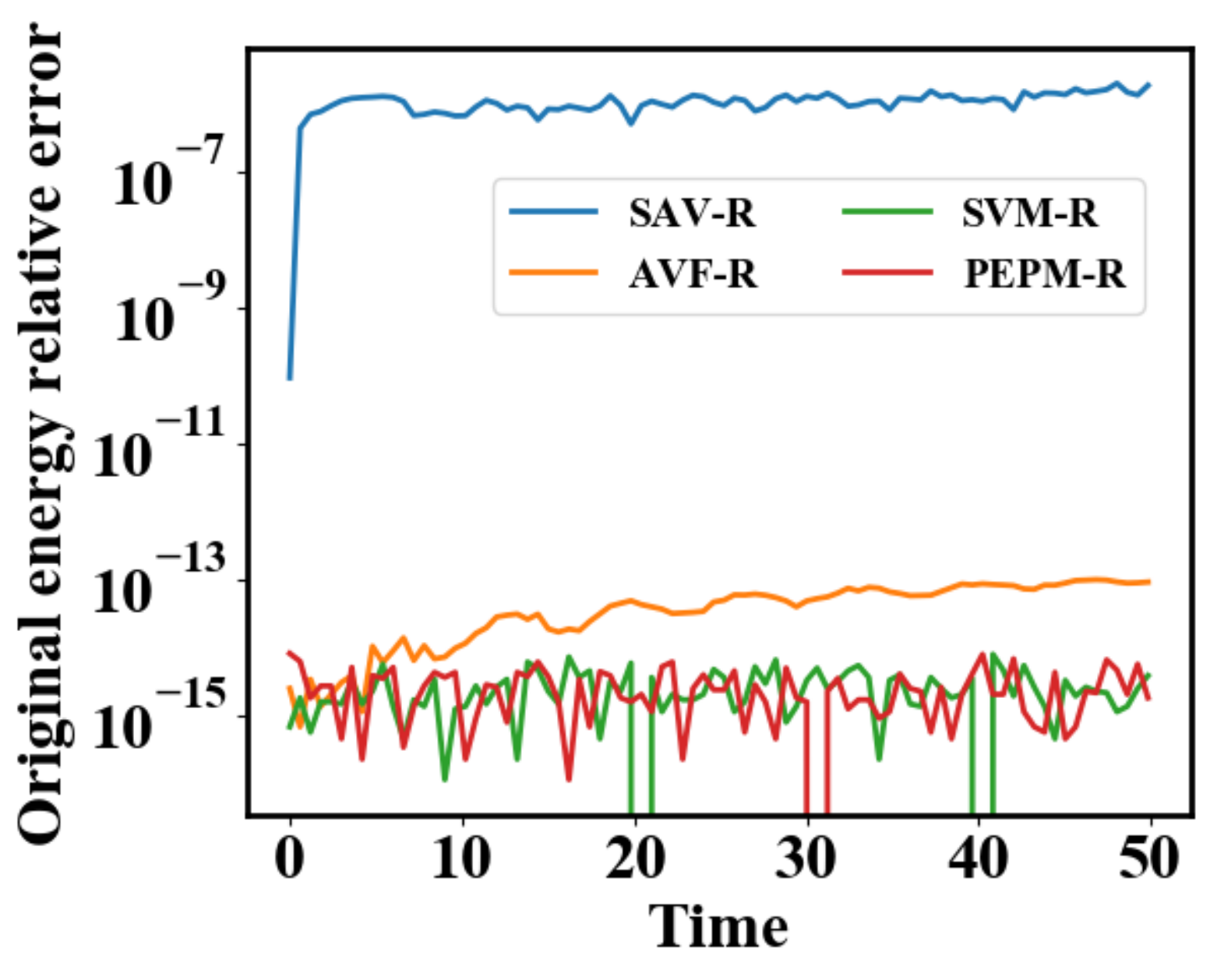}
		}
		\caption{{\bf Example \ref{eg:PLS}}:  Time evolution of  energy errors  using different schemes with $\tau = 0.01$ and $N_x = N_y =128$. The curves of the energy errors indicate our proposed schemes present a better conservation of the original energy as well as  the traditional energy-preserving algorithms i.e., {\bf AVF-M/-R}, but fails to {\bf SAV-M/-R} based on SAV approach.
			\label{fig:different-energy-error}}
	\end{figure}
\end{example}

\begin{example}[Line soliton in an inhomogeneous medium]\label{eg:LSIM}
	In this test, we consider  an inhomogeneous on large-area Josephson junction  given by the Josephson current density
	\begin{align*}
	\phi(x,y) = 1 + \sech^2\left(\sqrt{x^2 + y^2}\right),
	\end{align*}
	and the initial conditions
	\begin{align}
	\begin{cases}
	u(\bx, 0) = 4 \tan^{-1}\left[ \exp\left( \dfrac{x - 3.5}{0.954}  \right)\right],\\
	v(\bx, 0) = 0.629  \sech \left[ \exp\left( \dfrac{x - 3.5}{0.954}  \right)\right].
	\end{cases}
	\end{align}
	This model  is discretized spatially using the cosine pseudo-spectral method with $N_x = N_y = 128$  in a domain $\Omega = [-7, 7] \times [-7, 7]$.	The profiles of numerical solutions $\sin(u/2)$ with time step $\tau = 0.01$ are summarized in Figure \ref{fig:inhomogeneous-evolution}. Qualitatively, these numerical phenomenons are consistent to this simulations  in \cite{argyris1991finite,christiansen1981numerical,dehghan2008numerical}.  
	The conservative energy results of the four algorithms are listed in Figure \ref{fig:inhomogeneous-energy-error} (a). It can be observed that the numerical energy obtained by the new algorithms is well conserved. Figure \ref{fig:inhomogeneous-energy-error} (b) depicts the plots of the supplementary variable $\beta$ and the Lagrange multipliers $\lambda$. We found that $\beta$ and $\lambda$ are  up to ${\cal O}(\tau^2)$ and ${\cal O}(\tau^3)$.    These results strongly support our claim that our proposed schemes can be applied to predict accurate the motion of line soliton in an homogeneous medium. 
	\begin{figure}[H]
		\centering
		\subfigure{
			\includegraphics[width=0.22\textwidth,height=0.22\textwidth]{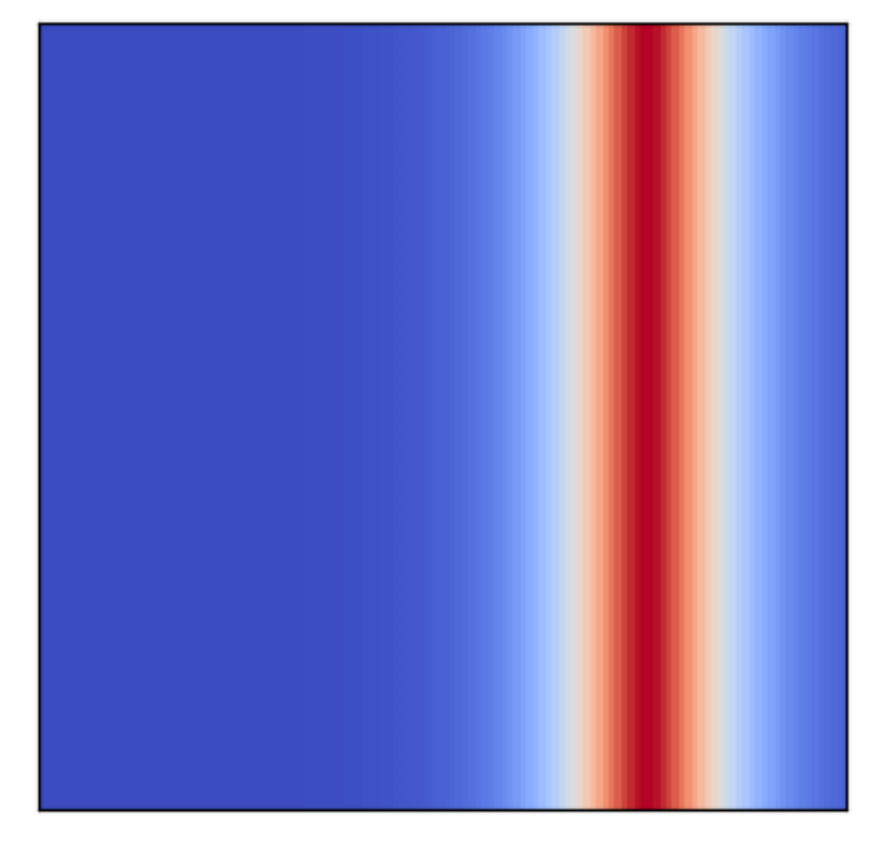}
			\includegraphics[width=0.22\textwidth,height=0.22\textwidth]{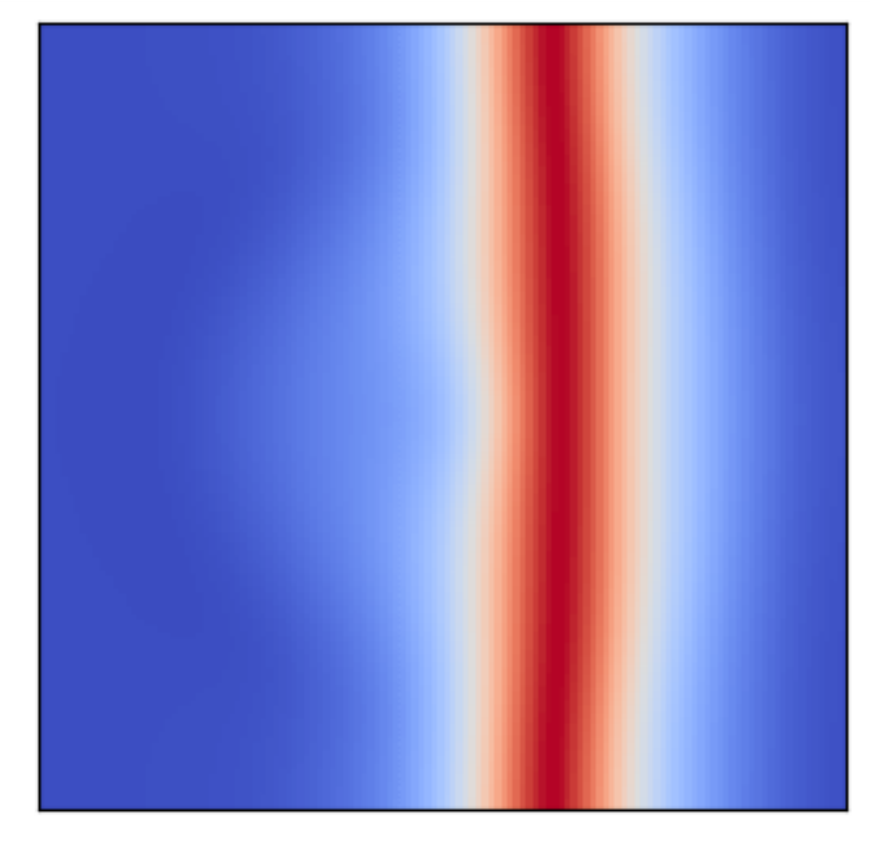}
			\includegraphics[width=0.22\textwidth,height=0.22\textwidth]{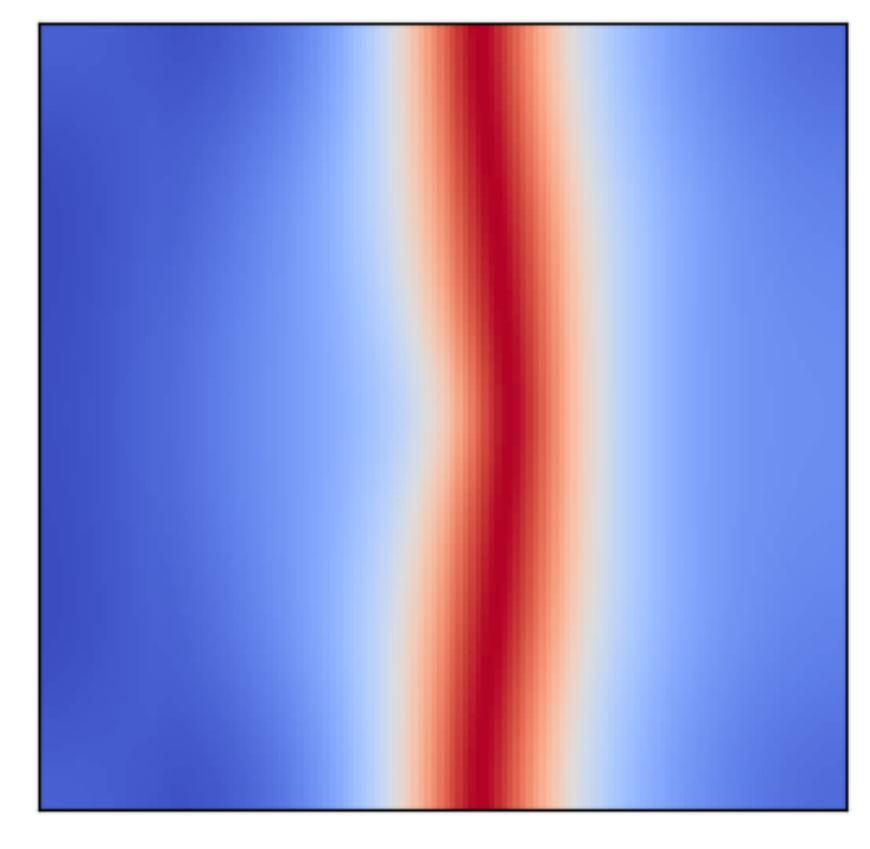}
			\includegraphics[width=0.22\textwidth,height=0.22\textwidth]{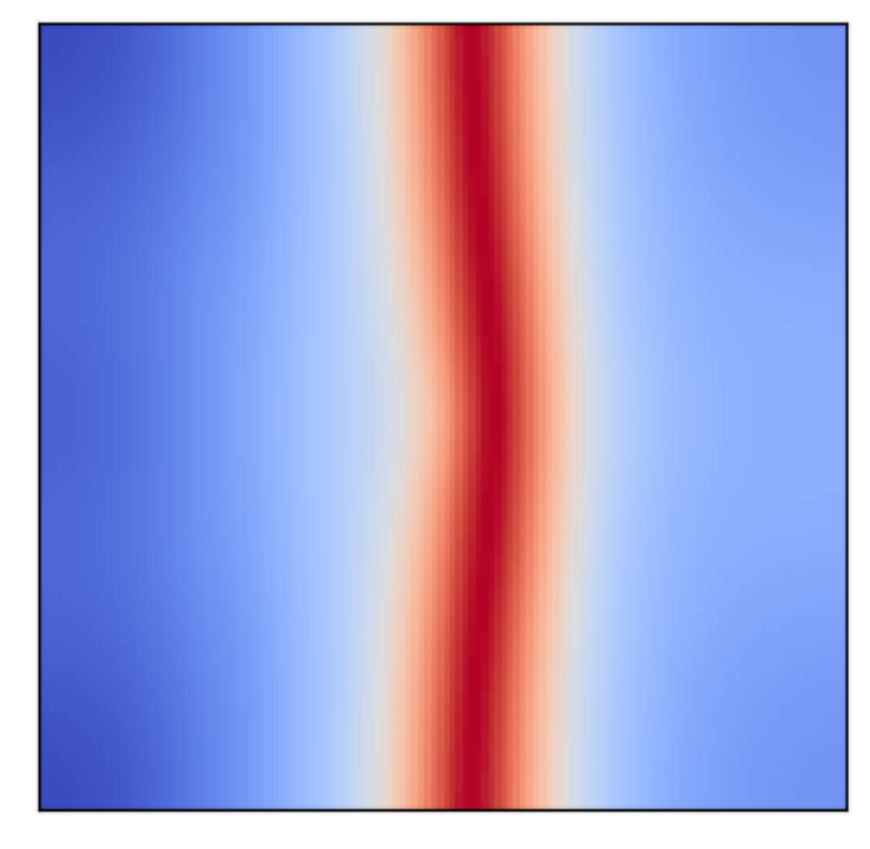}
		}
		\caption{{\bf Example \ref{eg:LSIM}}: The isolines of numerical solutions of $\sin(u/2)$. Snapshots are taken at $t = 0, 6, 12, 18$, respectively.   \label{fig:inhomogeneous-evolution}}
	\end{figure}
	
	\begin{figure}[H]
		\centering
		\subfigure{
			\includegraphics[width=0.38\textwidth,height=0.28\textwidth]{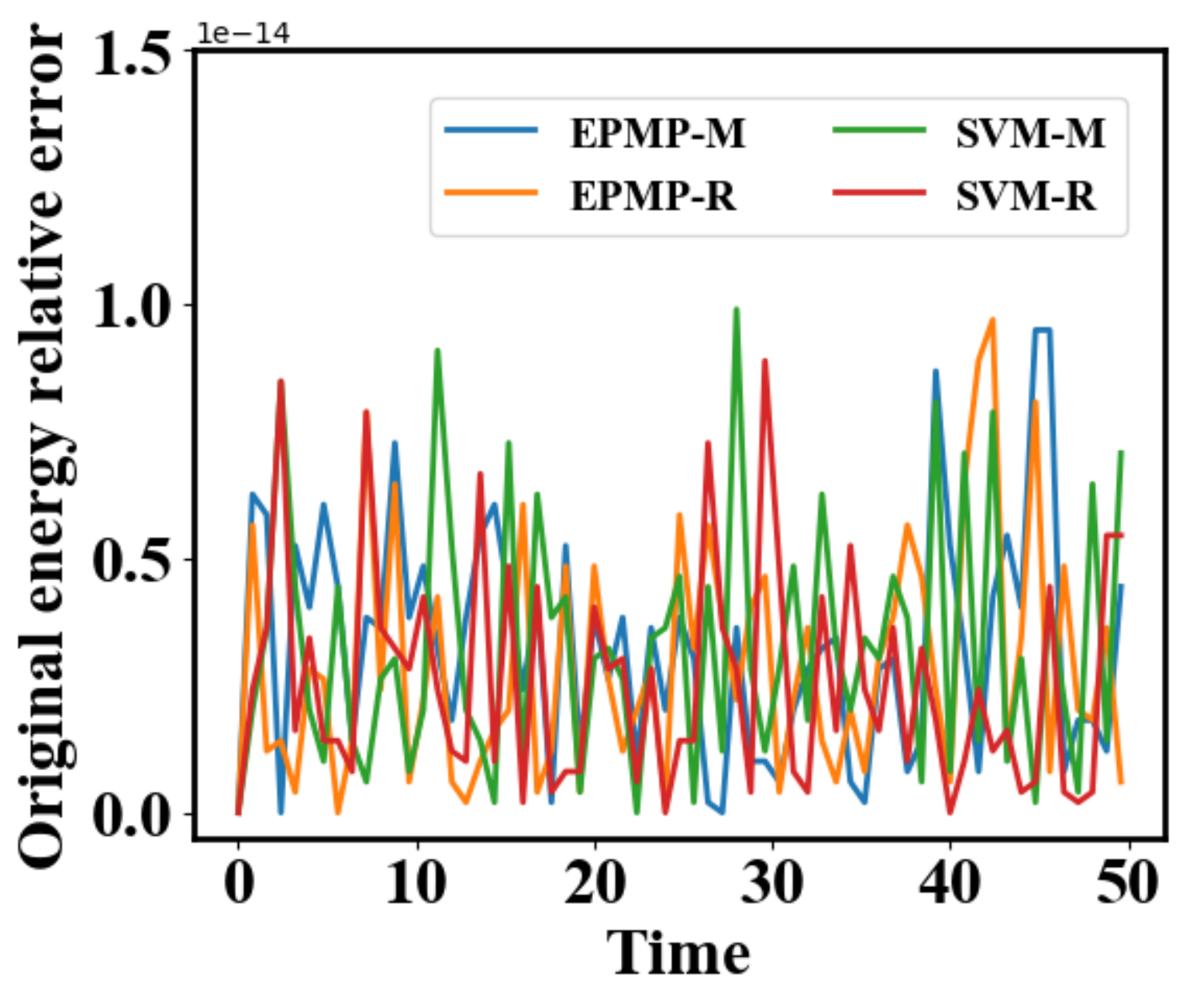}
		}
		\subfigure{
			\includegraphics[width=0.38\textwidth,height=0.28\textwidth]{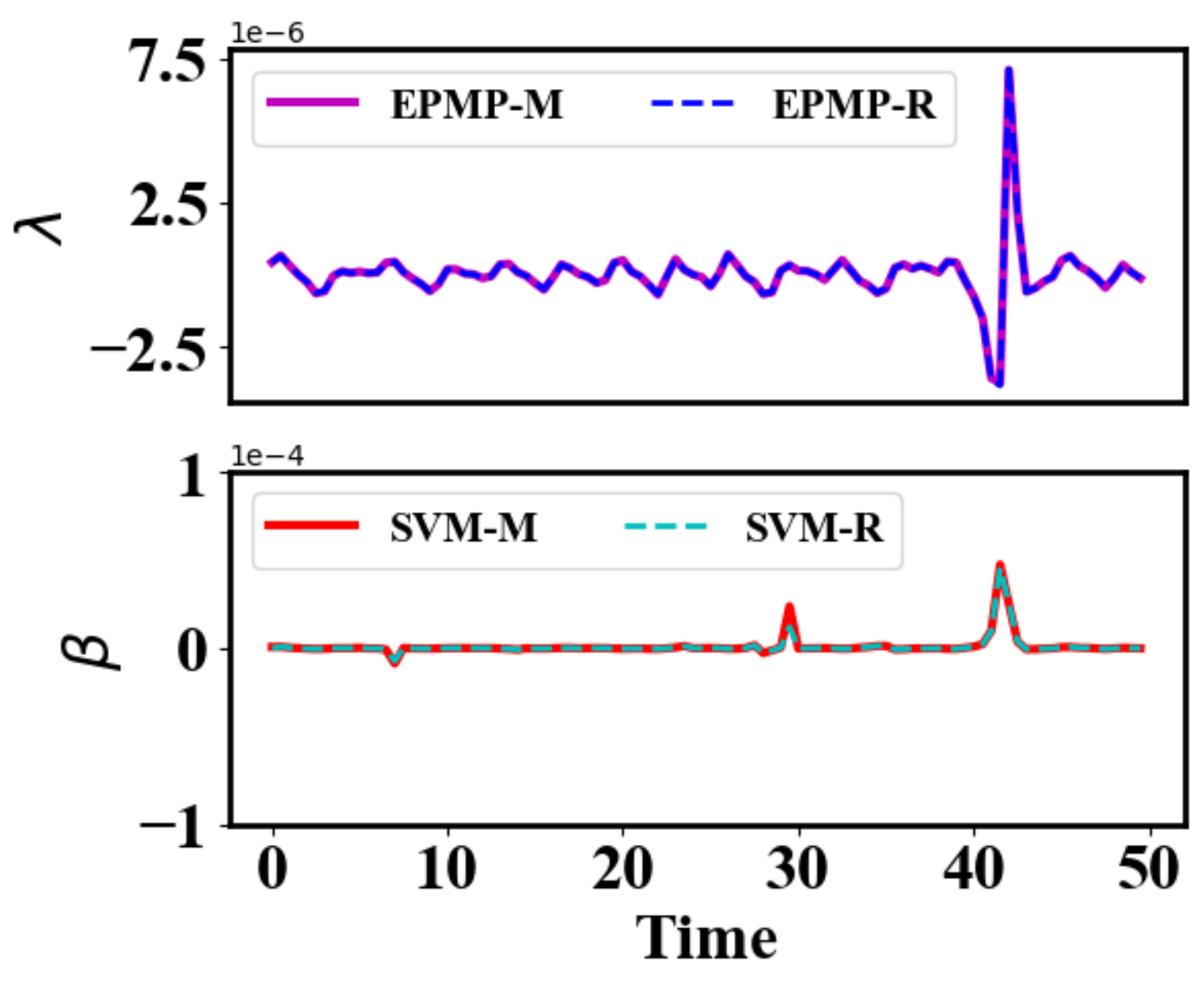}
		}
		\caption{{\bf Example \ref{eg:LSIM}}: Evolutions of the  energy errors using the four schemes with $\tau = 0.01$ and $N_x = N_y = 128$ (left). The curves of the energy errors show the four schemes obey the original energy conservation law.  Evolutions of the supplementary variable $\beta$ and the Lagrange multiplier $\lambda$ (right). \label{fig:inhomogeneous-energy-error}}
	\end{figure}
\end{example}
\subsection{Ring solitons}
\begin{example}[Circular ring soliton]\label{eg:CRS}
	In this test, we consider $\phi(x, y) = 1$ and  initial conditions are given by
	\begin{align*}
	\begin{cases}
	u(\bx, 0) = 4 \tan^{-1}\left[ \exp(3 - \sqrt{x^2 + y^2})  \right],\\
	v(\bx, 0) = 0.
	\end{cases}
	\end{align*}
	Following \cite{sheng2005numerical}, we choose the domain as $[-14, 14]^2$. We solve this SG model subject to Neumann boundary condition using cosine pseudo-spectral method with $128^2$ meshes. All numerical schemes, i.e.,  {\bf PEPM-M/-R} and {\bf SVM-M/-R},  are implemented. The corresponding numerical solution is presented in terms of $\sin(u/2)$ in Figure \ref{fig:circular}.  At the initial time, it can be seen that the ring soliton shrinks and as time goes on, oscillations and radiations begin to form and continue. These numerical phenomenons are consistent with the reported literatures \cite{sheng2005numerical,christiansen1981numerical}. The original energy errors are also plotted in Figure \ref{fig:circular_energy_error} (a), which indicates that the proposed schemes warrant the discrete energy  to round-off errors. In Figure \ref{fig:circular_energy_error}, we plot evolution of the supplementary variable $\beta(t)$  and the Lagrange multiplier $\lambda(t)$. We found that the maximum values of $\beta$ and $\lambda$ is up to $10^{-4}$ and $10^{-6}$ and they remain close to zero.
	In a word, the above numerical behaviors support our claim that our proposed schemes are very efficient to deal with the Neumann boundary conditions.
	
	\begin{figure}[H]
		\centering
		\subfigure{
			\includegraphics[width=0.3\textwidth,height=0.25\textwidth]{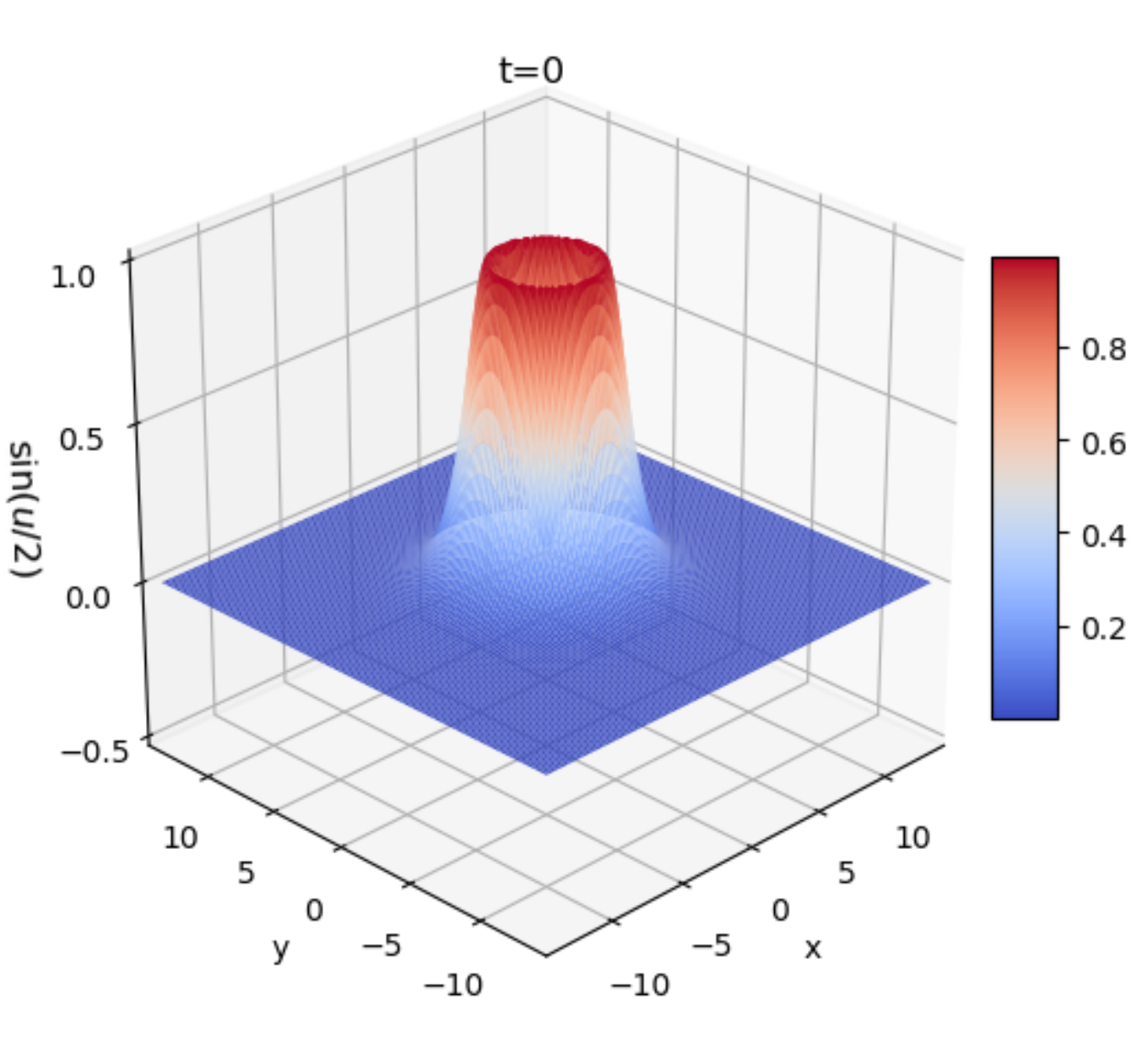}
		}
		\subfigure{
			\includegraphics[width=0.3\textwidth,height=0.25\textwidth]{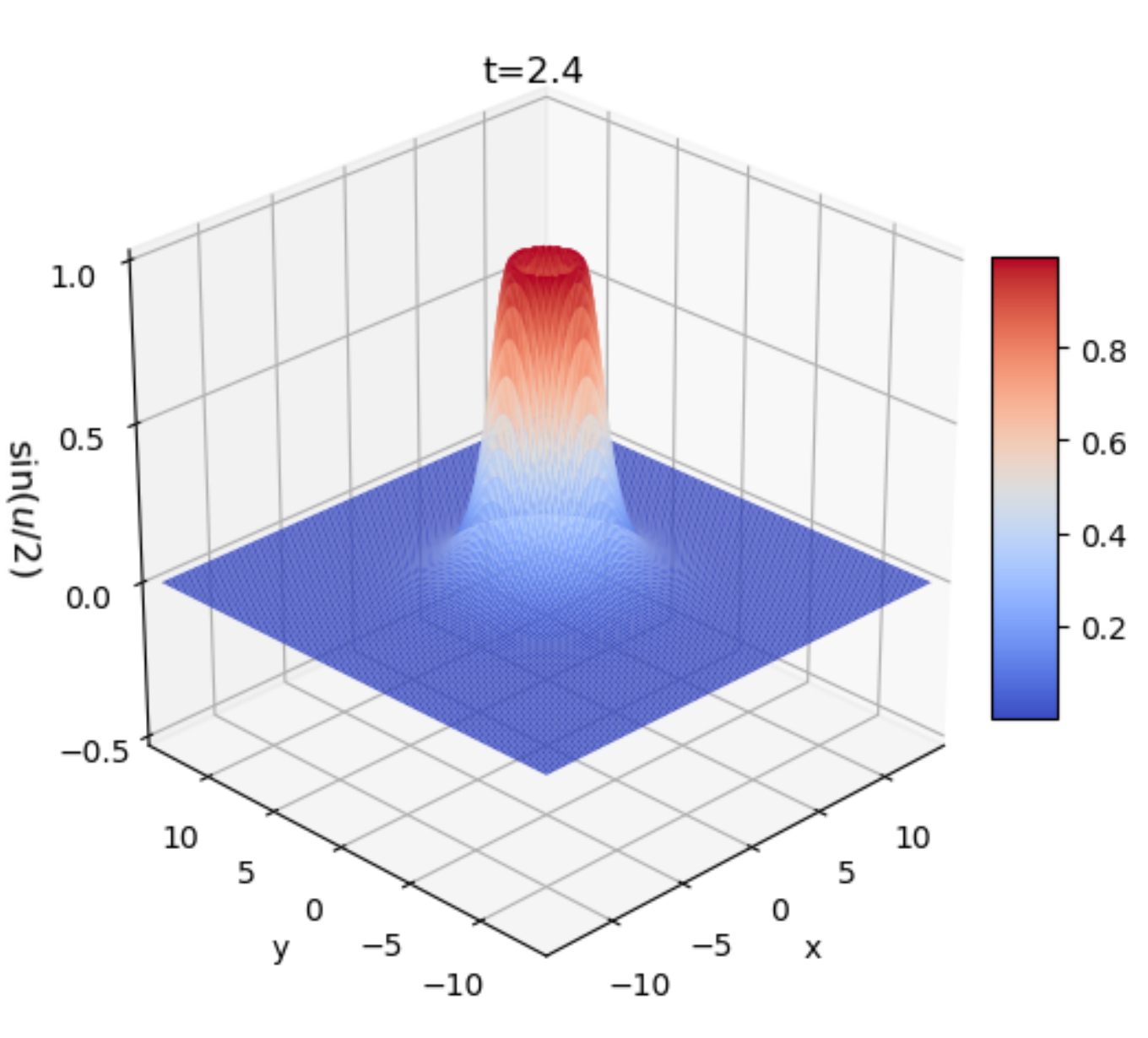}
		}
		\subfigure{
			\includegraphics[width=0.3\textwidth,height=0.25\textwidth]{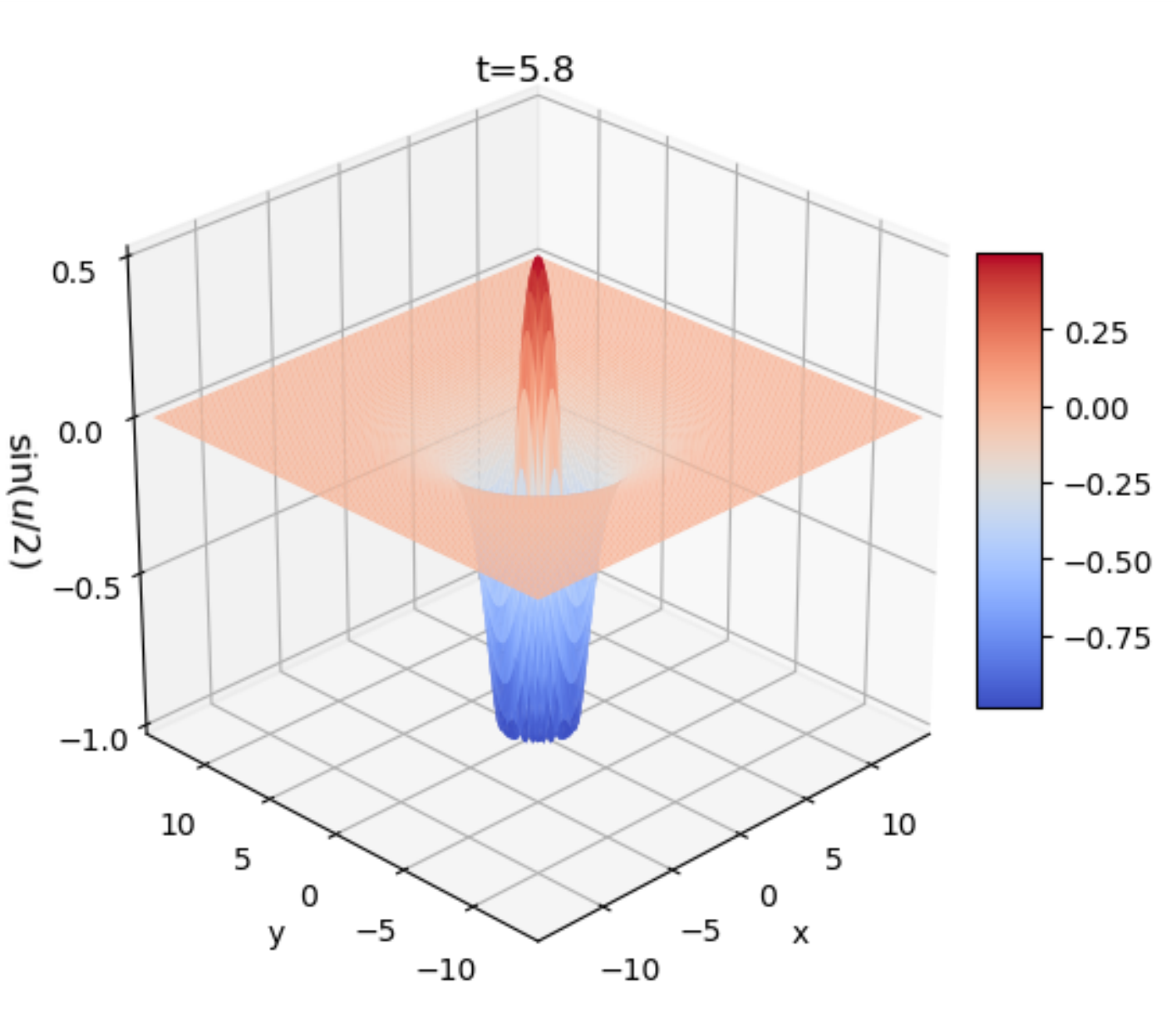}
		}
		\subfigure{
			\includegraphics[width=0.3\textwidth,height=0.25\textwidth]{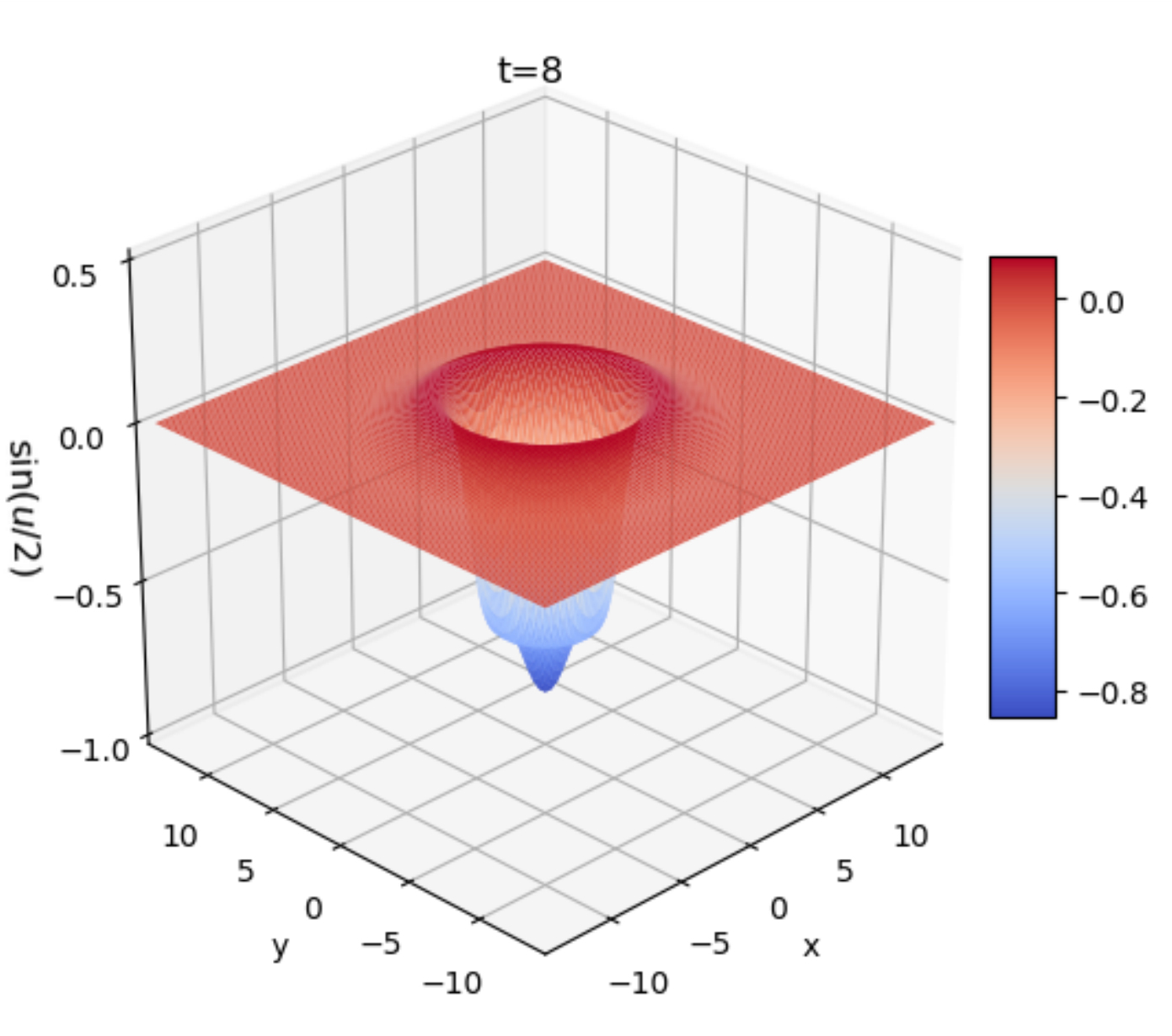}
		}
		\subfigure{
			\includegraphics[width=0.3\textwidth,height=0.25\textwidth]{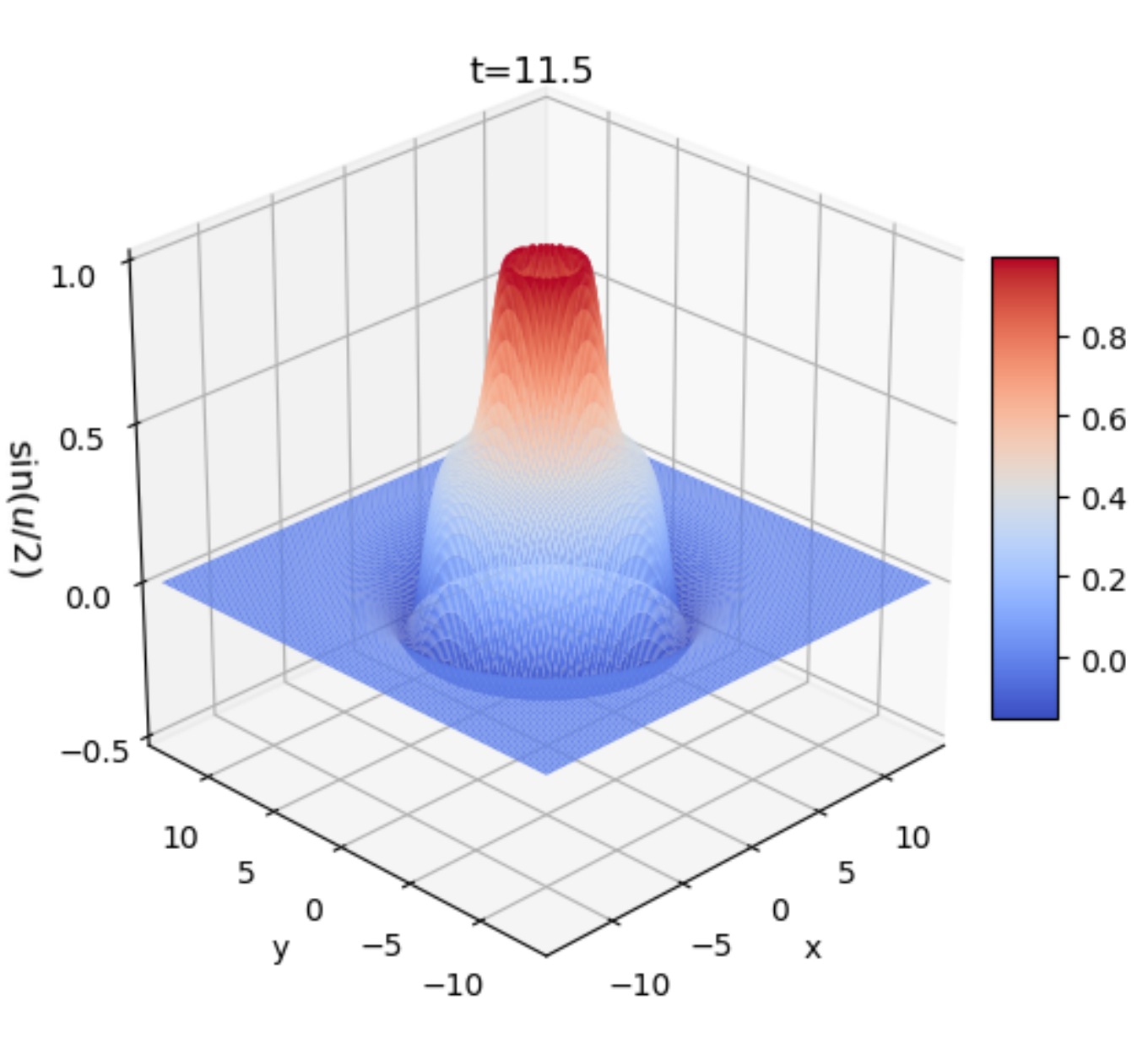}
		}
		\subfigure{
			\includegraphics[width=0.3\textwidth,height=0.25\textwidth]{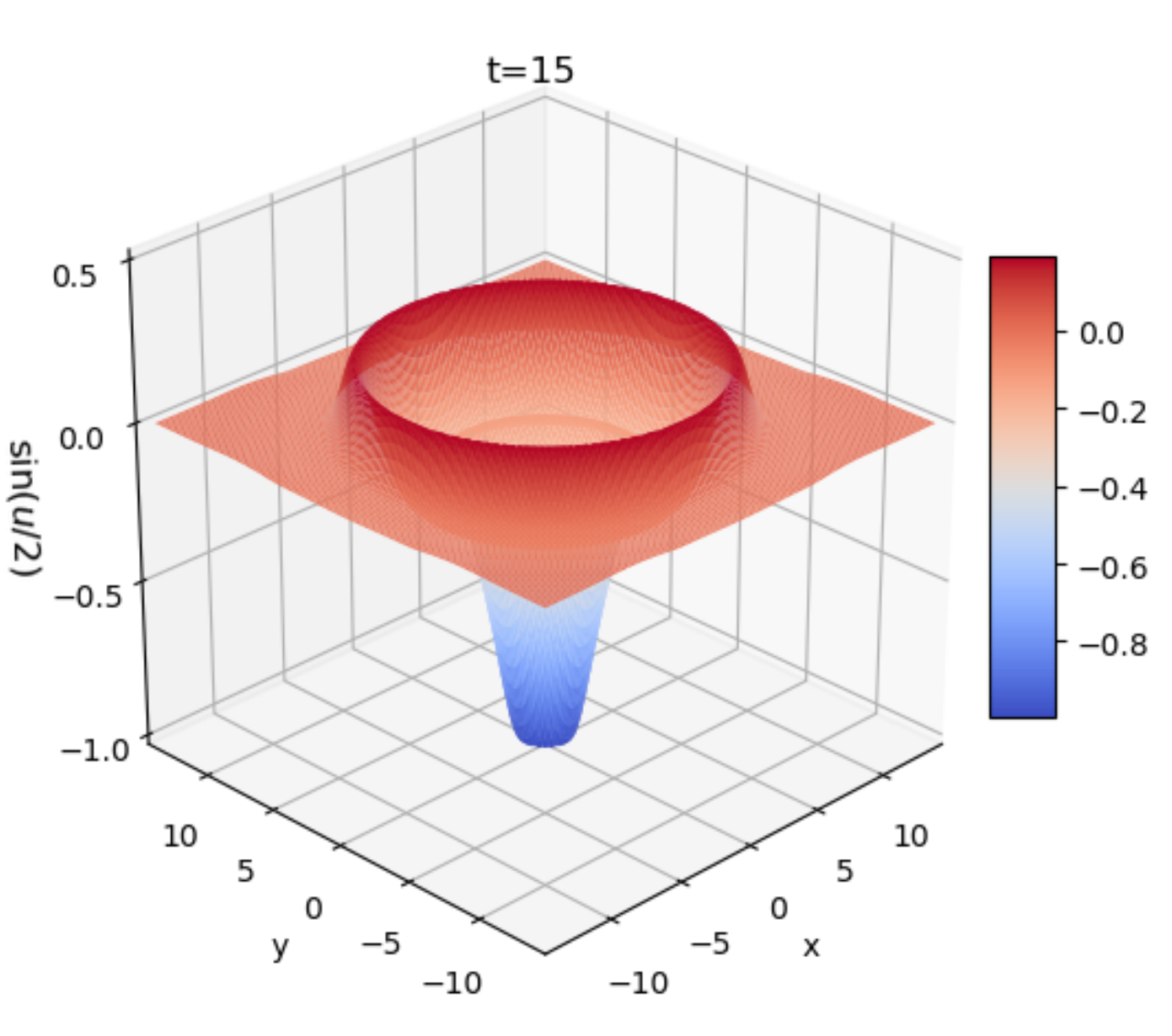}
		}
		\caption{{\bf Example \ref{eg:CRS}}: Evolutions of the  circular ring solitons in terms of $\sin(u/2)$ with $\tau = 0.01$ and $N_x = N_y =128$. Snapshots are taken at $t = 0, 2.4, 5.8, 8, 11.5, 15$, respectively.   \label{fig:circular}}
	\end{figure}
	\begin{figure}[H]
		\centering
		\subfigure{
			\includegraphics[width=0.40\textwidth,height=0.30\textwidth]{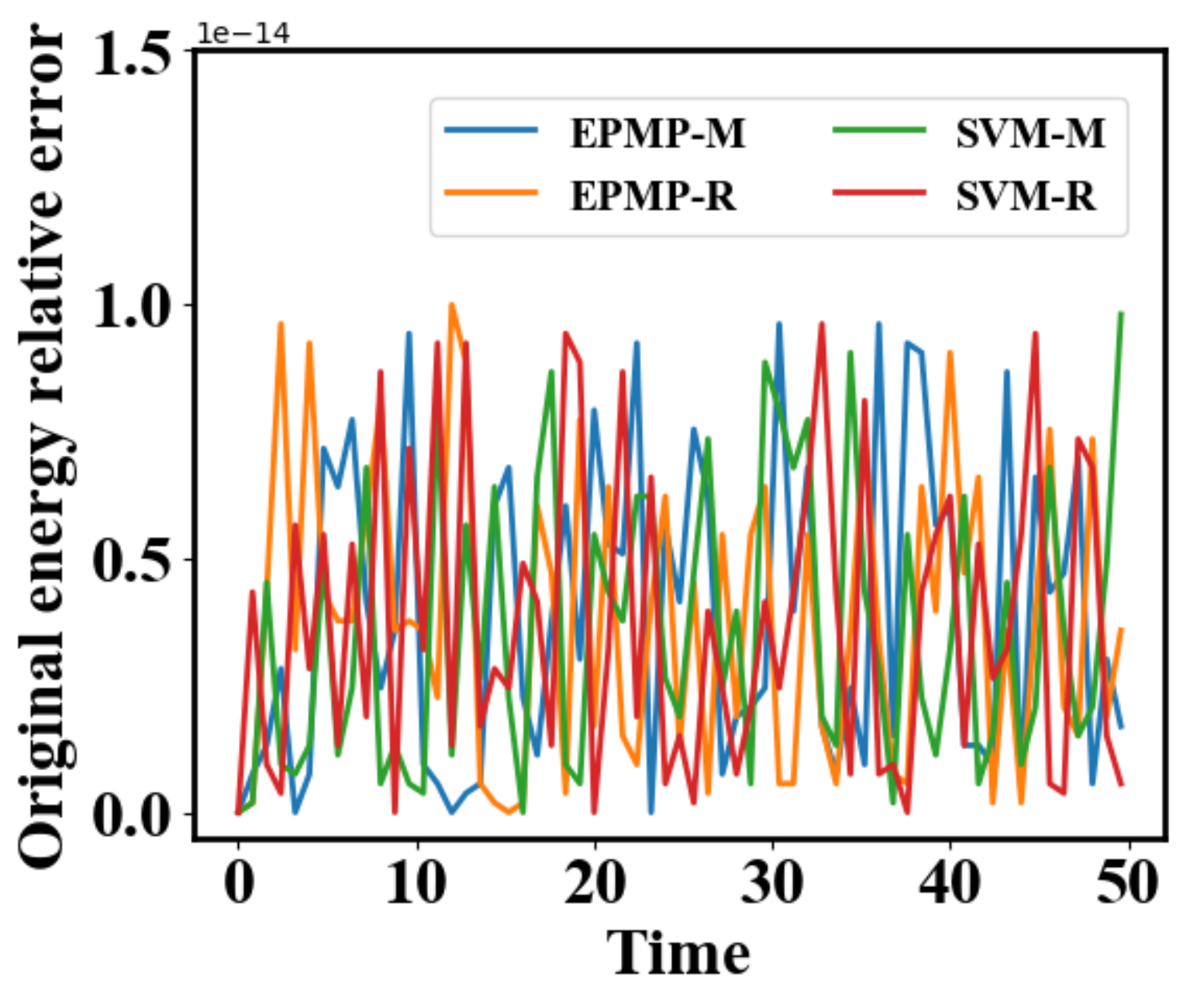}
		}
		\subfigure{
			\includegraphics[width=0.40\textwidth,height=0.30\textwidth]{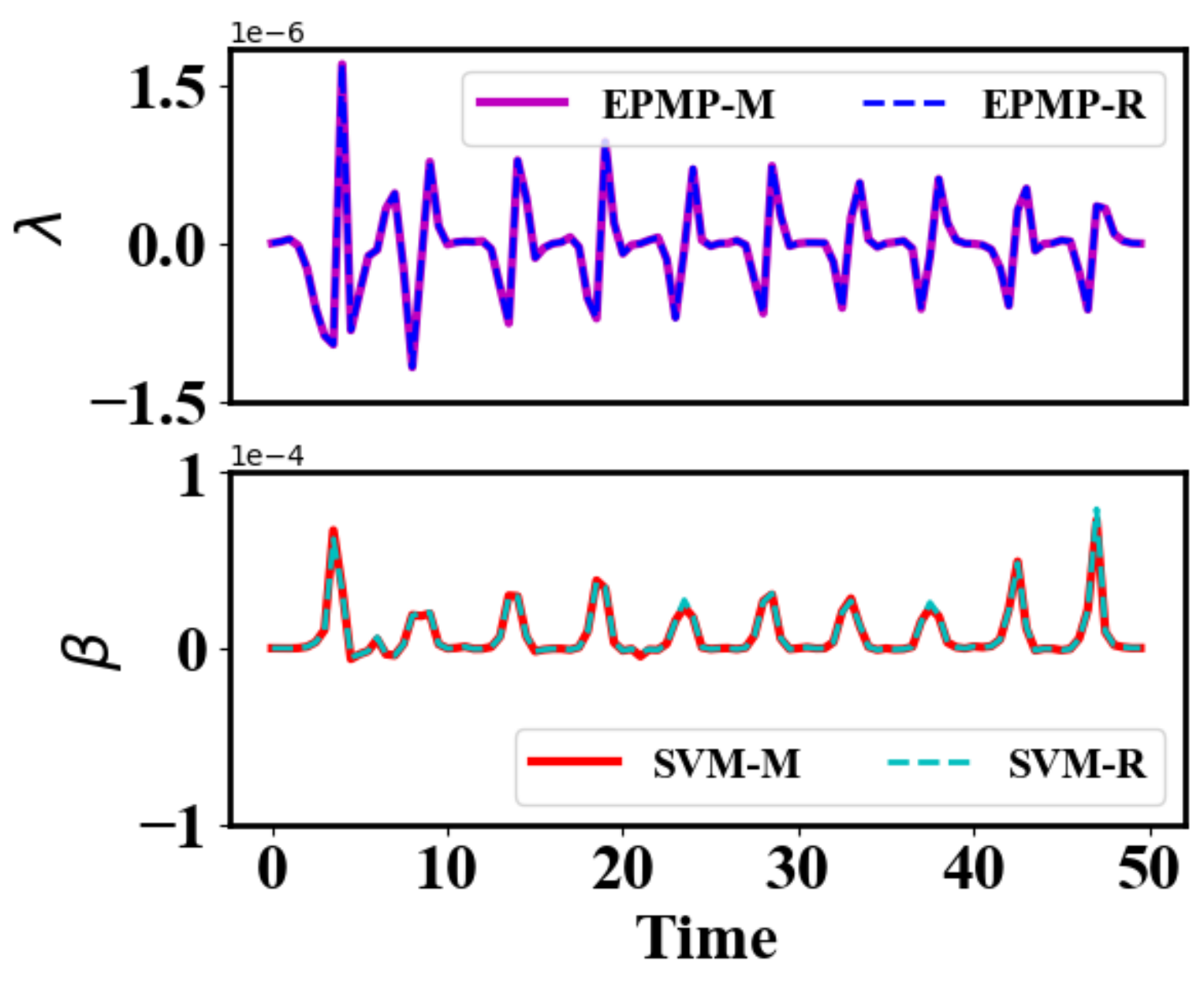}
		}
		\caption{{\bf Example \ref{eg:CRS}}: Evolutions of the  energy errors using the four schemes with $\tau = 0.01$ and $N_x = N_y =128$. The curves of the energy errors show the four schemes preserve the original energy conservation law very well (left).  Evolutions of the supplementary variable $\beta$ and the Lagrange multiplier $\lambda$. This subfigure shows $\beta(t)$ and $\lambda(t)$ may fluctuate near zero (right). \label{fig:circular_energy_error}}
	\end{figure}
\end{example}

\begin{example}[Collision of four circular solitons]\label{eg:CFCS}
	Finally, we end up with collisions of four expanding circular ring solitons, we select $\phi(x,y) = 1$ and initial conditions
	\begin{align*}
	& u(\bx, 0) = 4 \tan^{-1}\left[ \exp\left(\frac{4 - \sqrt{(x + 3)^2 + (y + 3)^2}}{0.436}\right)  \right],\\
	& v(\bx, 0) = 4.13\sech\left[ \exp\left(\frac{4 - \sqrt{(x + 3)^2 + (y + 3)^2}}{0.436}\right)  \right].
	\end{align*}
	This simulation  is based on an extension across $x = -10$ and $y = -10$ due to the symmetry. The computational domain is $[-30, 10] \times [-30, 10]$ and $128 \times 128$ grid points are used to discretize the space. 
	In this example, we intend to investigate that how affect the numerical behaviours via choosing the different supplementary functions $g[u, v]$. To save space, we just take the scheme {\bf SVM-M} as an example to demo, where $g_1[u, v] = \phi \sin u$ and $g_2[u, v] = \Delta u - \phi \sin u$.
	Figure \ref{fig:four-circular} shows the profiles of  $sin(u/2)$ by using {\bf SVM-M}, which demonstrates the collision between four expanding circular ring solitons in which the smaller ring solitons bounding an annular region emerge into a large ring soliton. The numerical behaviors agree qualitatively well with those in published literatures, for instance \cite{sheng2005numerical,argyris1991finite}. It is clear that the scheme simulate this problem very well. The changes of energy and the supplementary variable are list in Figure
	\ref{fig:four_circular_energy_error}. We observe that {\bf SVM-M} with two supplementary functions $g_1$ and $g_2$
	can warrant the original energy conservation law very accurately. 
	In addition, the supplementary variable $\beta(t)$  remains close to zero. 
	
	\begin{figure}[H]
		\centering
		\subfigure{
			\includegraphics[width=0.30\textwidth,height=0.25\textwidth]{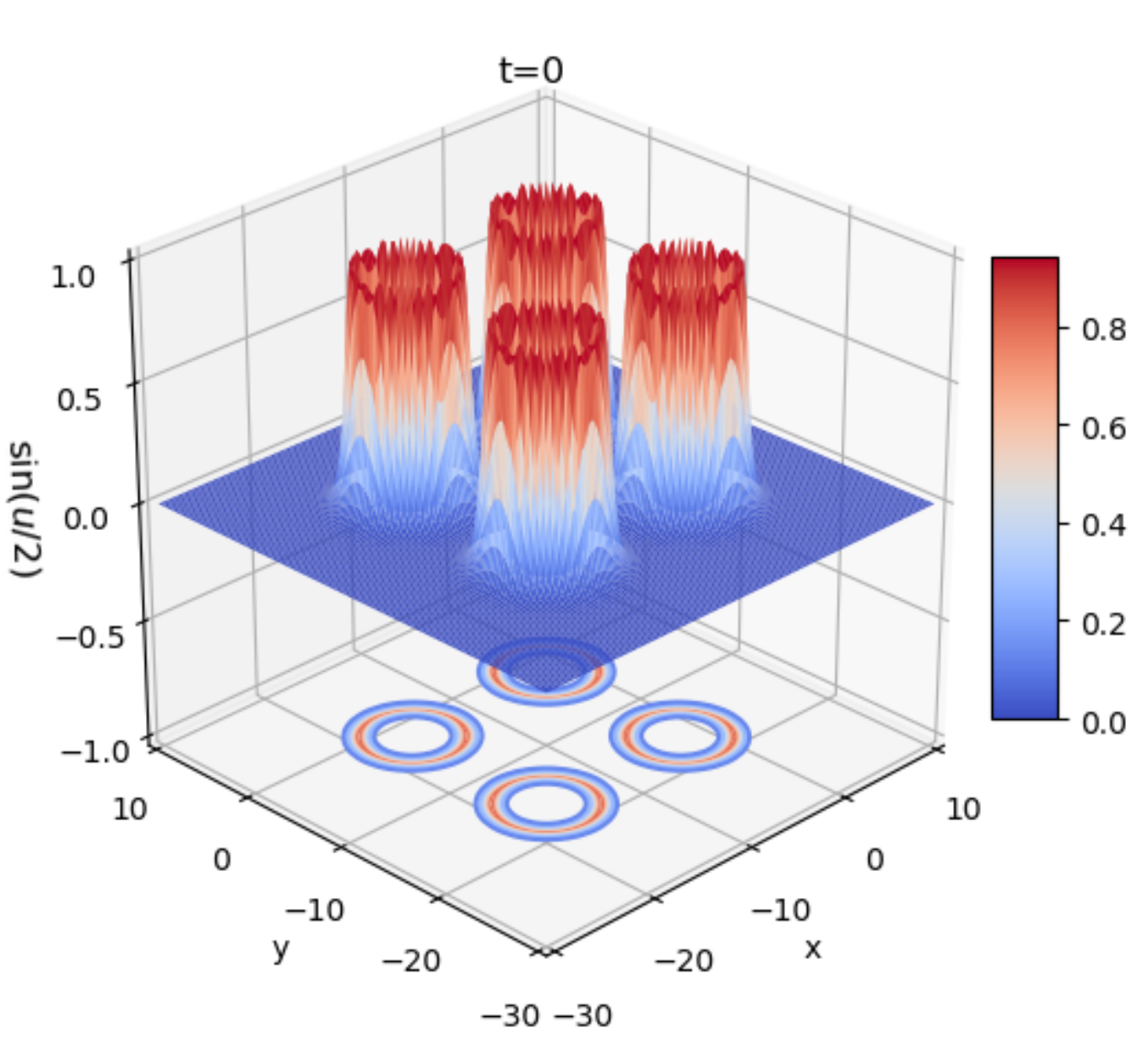}
			\includegraphics[width=0.30\textwidth,height=0.25\textwidth]{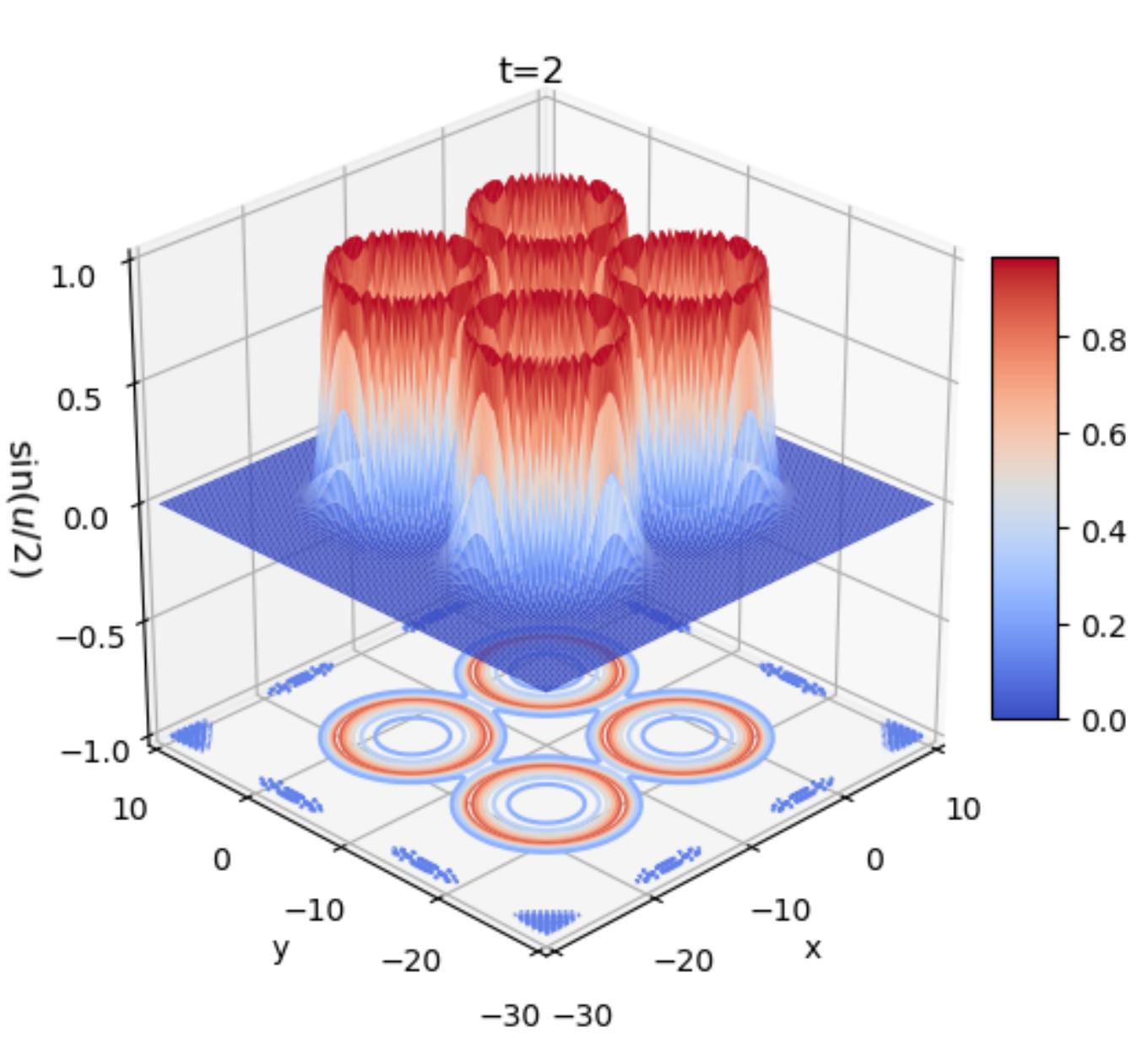}
			\includegraphics[width=0.30\textwidth,height=0.25\textwidth]{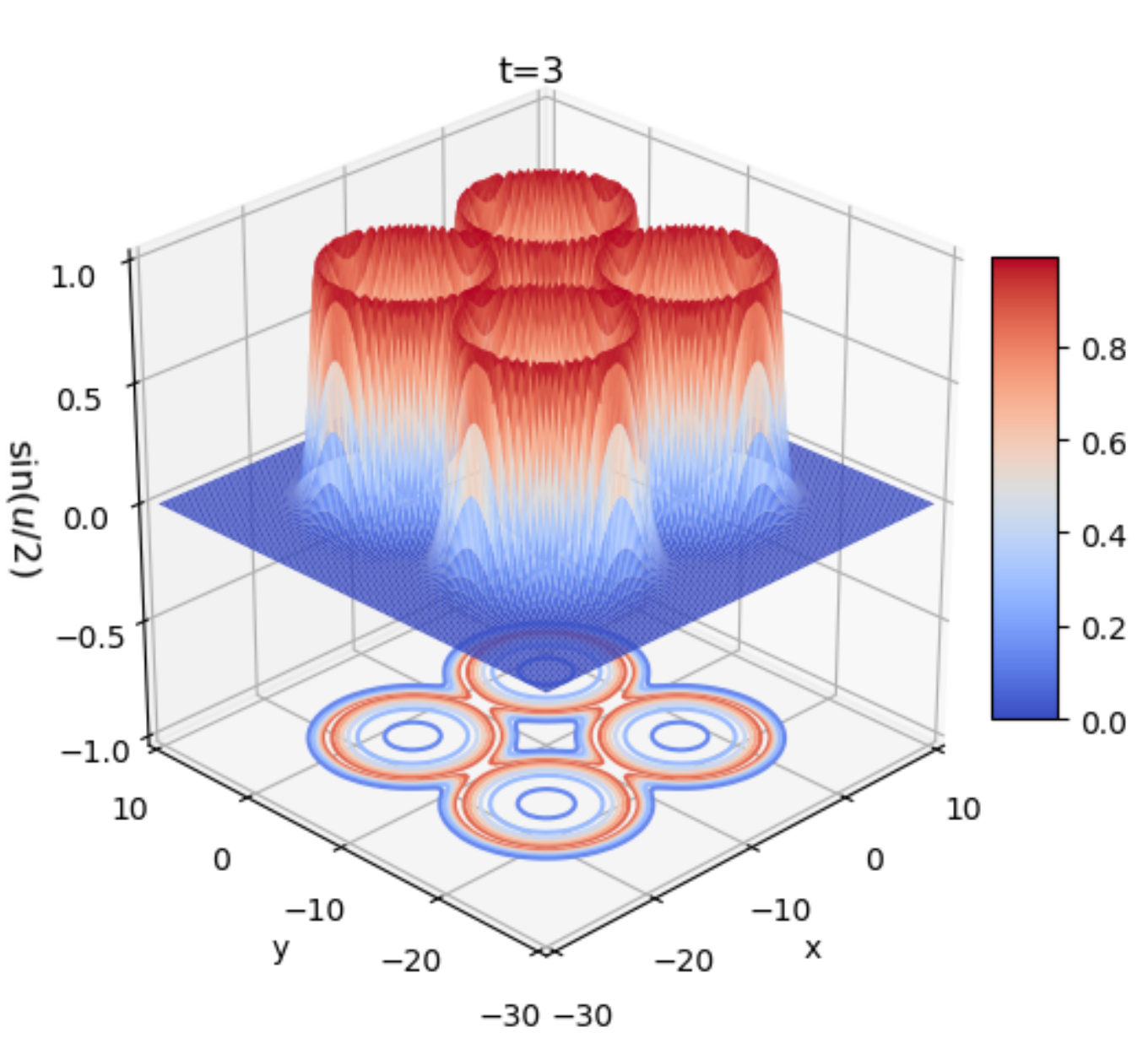}
		}
		\subfigure{
			\includegraphics[width=0.30\textwidth,height=0.25\textwidth]{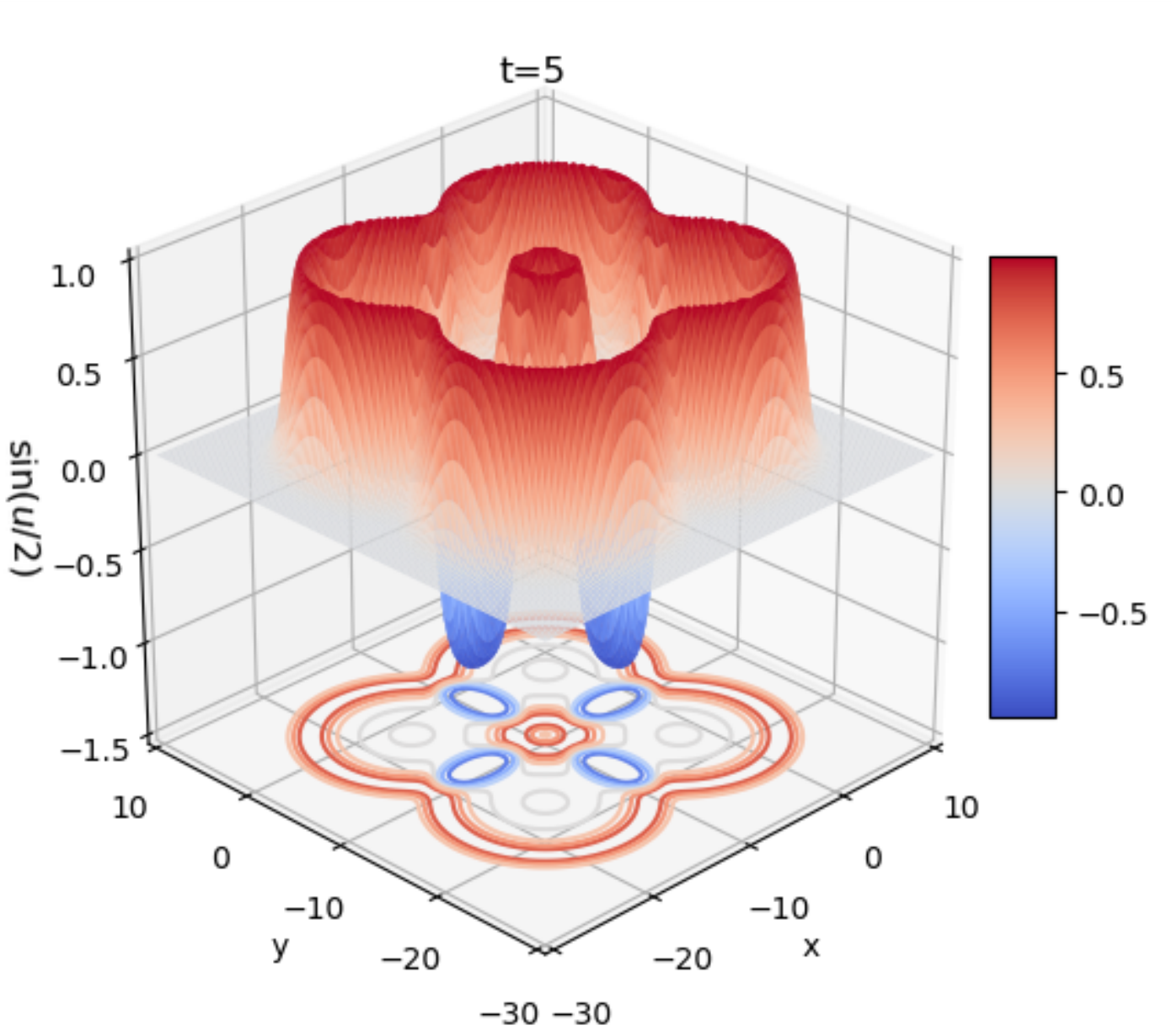}
			\includegraphics[width=0.30\textwidth,height=0.25\textwidth]{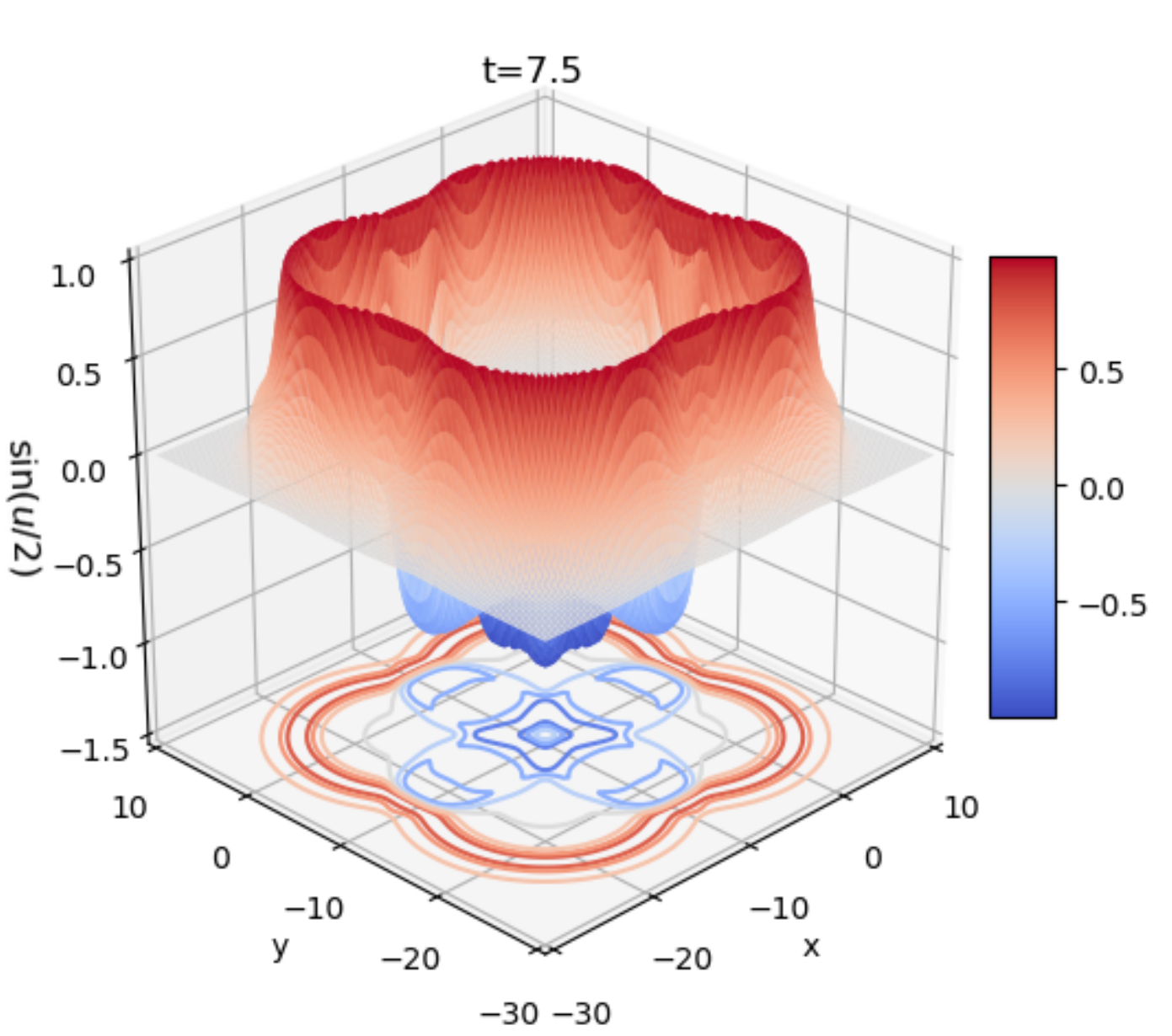}
			\includegraphics[width=0.30\textwidth,height=0.25\textwidth]{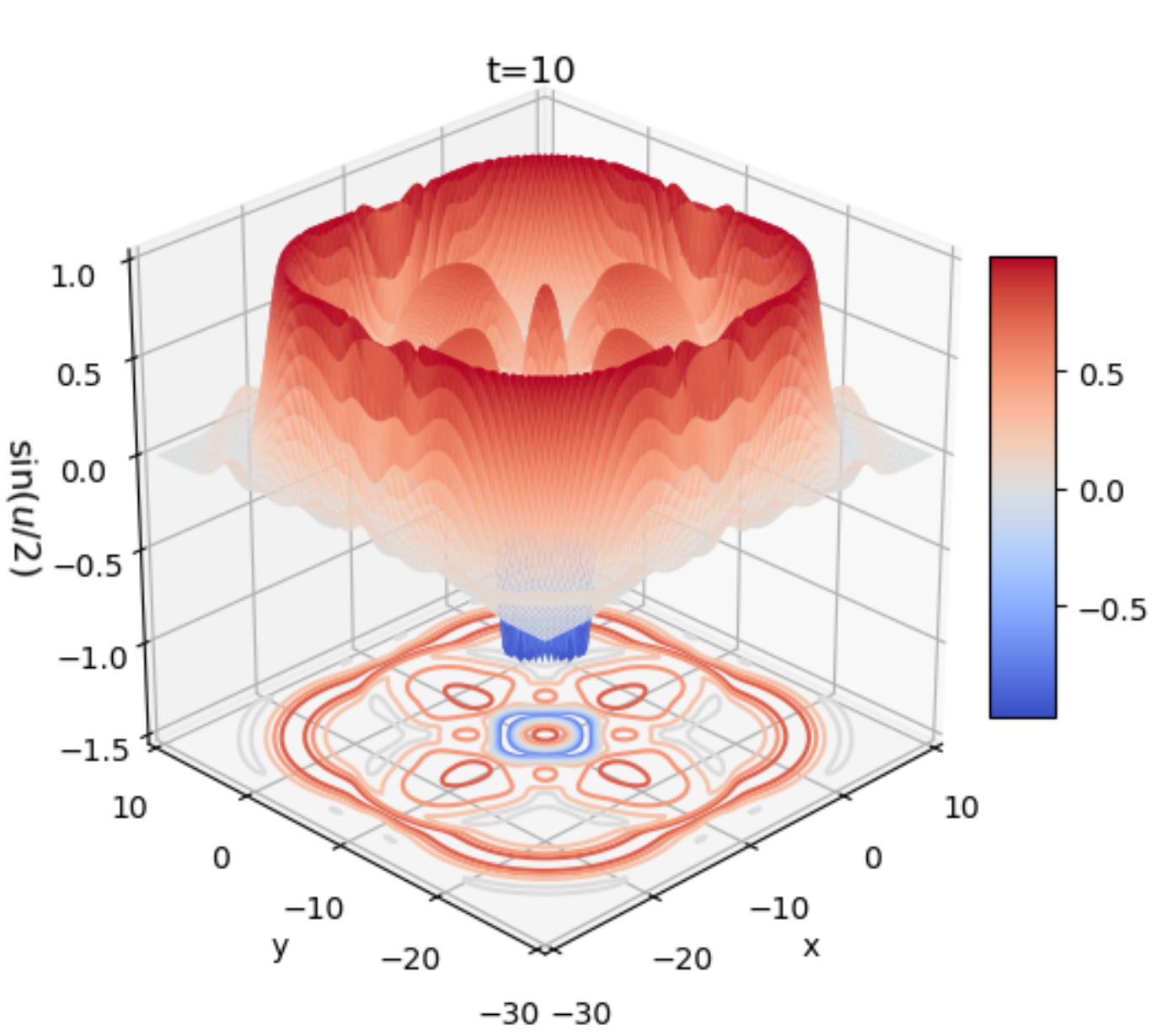}
		}
		\caption{{\bf Example \ref{eg:CFCS}}: Evolutions of the  circular ring solitons in terms of $\sin(u/2)$ using {\bf PEPM-M}  with $\tau = 0.01$ and $N_x = N_y =128$. Snapshots are taken at $t = 0, 2, 3, 5, 7.5, 10$, respectively.   \label{fig:four-circular}}
	\end{figure}

	\begin{figure}[H]
		\centering
		\subfigure{
			\includegraphics[width=0.40\textwidth,height=0.30\textwidth]{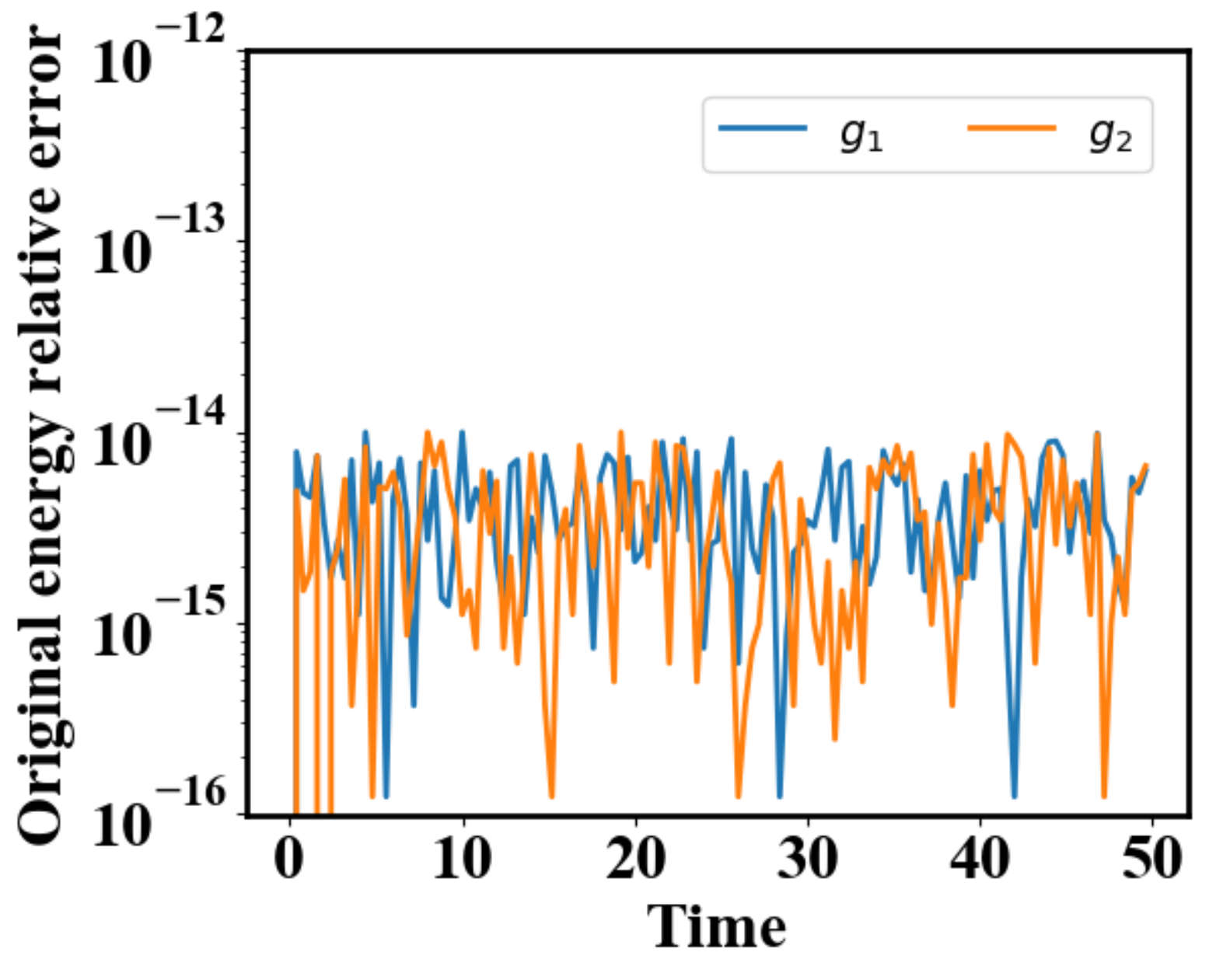}
		}
		\subfigure{
			\includegraphics[width=0.40\textwidth,height=0.30\textwidth]{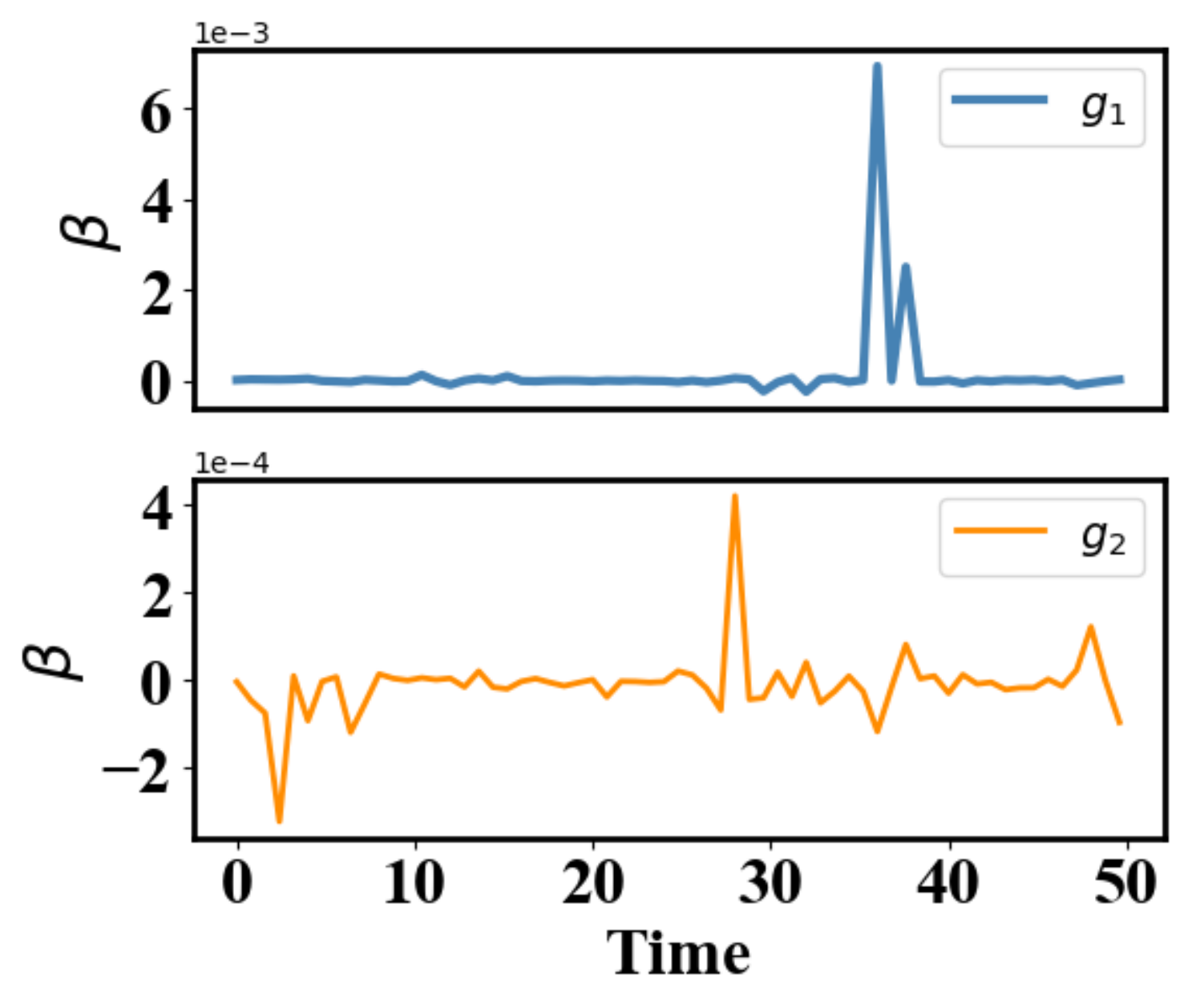}
		}
		\caption{{\bf Example \ref{eg:CFCS}}: Evolutions of the  energy errors using different  function $g_1$ and $g_2$ with $\tau = 0.01$ and $N_x = N_y =128$. The curves of the energy errors show {\bf SVM-M} with  two different supplied functions preserves the original energy conservation law very well (left).   This subfigure of  the supplementary variable shows $\beta(t)$ is about ${\cal O}(\tau^2)$ (right). \label{fig:four_circular_energy_error}}
	\end{figure}

\end{example}

\section{Conclusions}\label{sec:conclusions}
We have developed two classes of newly structure-preserving algorithms for sine-Gordon  equation subject to Neumann boundary conditions, which are based on the projection approach and the supplementary variable method.  Two kinds of  cosine pseudo-spectral  fashions are proposed to achieve high-order in the spatial discretization, one is based on the mid-point grid and the other is on the regular grid, which admit that the semi-discrete system both are canonical Hamiltonian structure and energy conservation law. In the first approach, we present a prediction-correction Crank-Nicolson method based on an energy projection skill to derive efficient fully-discrete schemes. Although they retain second-order in time, the error is much smaller than other energy-preserving algorithms, and the 2nd order convergence rate in time and high-order accuracy in space are attained. In the second approach, we proposed a novel idea utilizing a SVM to develop energy-preserving algorithms at two different grids for the SG model subject to Neumann boundary conditions. The resulting numerical schemes are based on the blend of the cosine pseudo-spectral method in space and the linearized Crank-Nicolson method  in time. Moreover, these schemes are shown to possess the discrete energy conservation law.
In addition, the proposed numerical schemes require to solve a scalar nonlinear equation by a Newton iteration, which is negligible  than the main computation cost. But they lead to additional difficulty in its convergence and error analysis, which is our ongoing project. Numerical tests with benchmark problems are shown to illustrate the accuracy and effectiveness of the proposed schemes. The idea and methodology developed here can be extended to more general Hamiltonian systems and constructed high-order algorithms, which will be reported in a sequel.

\section*{Acknowledgements} 
The research is partially supported by the China Postdoctoral Science Foundation through Grant 2020M670116,
the Foundation of Jiangsu Key Laboratory for Numerical Simulation of Large Scale Complex Systems (202001, 202002), the Natural Science Foundation of Jiangsu Province (Grant No. BK20180413),  the National Key Research and Development Project of China (Grant No. 2016YFC0600310, 2018YFC0603500, 2018YFC1504205) and the National Natural Science Foundation of China (Grant No. 11771213,  11801269 and  NSAF-U1930402).

\end{document}